\documentclass[a4paper,11pt]{article}
\usepackage[centertags]{amsmath}
\usepackage{amssymb}
\usepackage{amsthm}
\usepackage{graphicx}
\usepackage{mathdots}
\usepackage{amsfonts}
\usepackage{mathrsfs}
\usepackage{times,cite}
\usepackage{geometry}
\usepackage{color}
\newtheorem{thm}{Theorem}[section]

\newtheorem{prop}[thm]{Proposition}

\theoremstyle{definition}
\newtheorem{rem}[thm]{Remark}

\newenvironment {Proof} {\noindent {\bf Proof. }}{\hfill $\square$\par\vspace{3mm}}
\newenvironment {Proof of Theorem} {\noindent {\bf Proof of Theorem 2.1}}{\hfill $\square$\par\vspace{3mm}}
\numberwithin{equation}{section}
\allowdisplaybreaks    

\newenvironment{acknowledgment}{\noindent{\bf Acknowledgment }}{}
\date{}
\title{Stability analysis for the incompressible Navier-Stokes equations with Navier boundary conditions}
\author{{ Shijin Ding$^a$,~ Quanrong Li$^b$\thanks{Corresponding author. },~ Zhouping Xin$^c$}\\
{\it\small $^a$School of Mathematical Sciences, South China Normal University,}\\
{\it\small Guangzhou, 510631, China}\\
{\it\small $^b$College of Mathematics and Statistics, Shenzhen University,}\\
{\it\small Shenzhen, 518060, China}\\
{\it\small $^c$The Institute of Mathematical Sciences and Department of Mathematics,}\\
{\it\small The Chinese University of Hong Kong, Shatin, N.T., Hong Kong.}\\
{\small E-mail addresses: {\it dingsj@scnu.edu.cn}(S. Ding), {\it quanrong\_li@szu.edu.cn}({Q. Li}),}\\
 \small{\it zpxin@ims.cuhk.edu.hk}(Z. Xin)}
\begin{document}
\newcommand{\D}{\displaystyle}
\maketitle

 \begin{abstract}
 This paper is concerned with the instability and stability of the trivial steady states of the incompressible Navier-Stokes equations with Navier-slip boundary conditions in a slab domain in dimension two. The main results show that the stability (or instability) of this constant equilibrium depends crucially on whether the boundaries dissipate energy and the strengthen of the viscosity and slip length. It is shown that in the case that when all the boundaries are dissipative, then nonlinear asymptotic stability holds true. Otherwise, there is a sharp critical viscosity, which distinguishes the linear and nonlinear stability from instability.
 \end{abstract}

{\bf Keywords: }stability and instability, Navier-Stokes equations, Navier boundary conditions,

 critical viscosity.

{\em AMS Subject Classification:} 76N10, 35Q30, 35R35.

\section{Formulation of the problem}

In this paper, we are interested in the following incompressible Navier-Stokes
equations
\begin{align}
\begin{cases}\label{0.1}
\partial_t{\mathbf v}+{\mathbf v}\cdot\nabla{\mathbf v}+\nabla p-\mu\Delta {\mathbf v}=0,\\
\mathrm{div}{\mathbf v}=0,
\end{cases}
\end{align}
in $\Omega, \ t>0$ where $\Omega$ is a domain in $R^N$ ($N\geq 2$), ${\mathbf v}$ is the velocity vector field, $p$ is the pressure.

System (\ref{0.1}) is mostly studied with no-slip boundary condition, i.e., Dirichlet boundary condition which means that the fluid does not slip along the boundary. However, this condition is not always realistic and leads to
induce a strong boundary layer in general. For example, hurricanes and tornadoes, do slip along the ground, lose energy as they slip and do not penetrate the ground (see \cite{3}). Other examples about the slip of the fluid on the boundary occur when moderate pressure is involved such as in high altitude aerodynamics (see \cite{haase}), or in immiscible two phase flows, the moving contact line is not compatible with no-slip boundary condition, see \cite{Amrouche} and \cite{guoyan2017}. As early as 1827, Navier \cite{18} had taken such factors into account and proposed a boundary condition as follows which is now called the Navier boundary condition in which there is a stagnant layer of fluid close to boundary allowing a fluid to slip
\begin{equation}\label{0.2}
{\mathbf v}\cdot{\mathbf n}=0 \ \ \mathrm{and}\ \ [(-p\mathbb{I}+\mu(\nabla {\mathbf v}+\nabla^T{\mathbf v}))\cdot {\mathbf n}]\cdot\tau=\alpha{\mathbf v}\cdot\tau,\ \mathrm{on}\ \partial\Omega\\
\end{equation}
where ${\mathbf n}$ is the outward normal vector field to $\partial\Omega$ and $\tau$ is the tangential vector. In Navier slip boundary condition (\ref{0.2}), $\alpha$ stands for
a physical meaning parameter which is either a constant or a function in $L^\infty(\partial\Omega)$ \cite{kelliher}, even a smooth matrix \cite{gie}. Here we restrict ourselves to the case that $\alpha$ is a constant

For such Navier bounadry value problems, the situation $\alpha\leq 0$ which reflects, in general, the friction between the fluid and the boundary, is the classical case and has got extensive attentions by physicists and mathematicians in studying the existence, uniqueness, regularity and vanishing viscosity to system (\ref{0.1}). The first pioneer paper on the mathematical rigorous analysis of the Navier-Stokes equation with Navier boundary conditions should be due to Solonnikov and \v{S}\v{c}adilov \cite{20} for the linearized stationary equations, while the existence of the weak solutions and regularity for the nonlinear case are obtained by B. da Veiga \cite{21} on half-space. Recently, \cite{2} and the references therein give some more specified results on existence and regularity of the solutions for various domains. In addition, for results on the vanishing viscosity limit for the evolutionary case, see \cite{23,24} and the references therein. For more physical applications and numerical analysis details, see \cite{1,4,6,13,14,17,19,191}.

Compared with the case $\alpha\leq 0$, the Navier boundary value problem with $\alpha>0$ has got less considerations. But just as Serrin stated in 1959 , $\alpha$ does not need to have defined sign (see Serrin \cite{191}, p.240). Of course, this is the case which does exist in reality and in physics. For example, for flat hybrid gas-liquid surfaces, the effective slip length $\alpha$ is always positive \cite{haase}, see also \cite{guoyan2017} for the contact line problem. Other examples are as follows. Navier boundary condition (\ref{0.2}) with $\alpha>0$ is also generally used for the simulations of flows in the presence of rough boundaries such as in aerodynamics, or in the case of permeable boundary where condition (\ref{0.2}) is called Beavers-Joseph's law (\cite{6}, \cite{Amrouche}), or in weather forecasts and in hemodynamics (\cite{6}, \cite{1993}), or when the boundary wall accelerates the fluid, see \cite{1995} and \cite{3}.

There are several kinds of concept on instability, we refer to \cite{7,8}. The most common concept is the Rayleigh-Taylor stability and instability due to heavier fluid on the upper forced by gravity, called RT stability and RT instability. RT stability and RT instability have been extensively studied, see \cite{9,10,11,12,15,16,22} and references therein. However, these researches for the stability problems are most subject to the no-slip boundary conditions. In 2016, Hailiang Li and Xingwei Zhang in \cite{lizhang} obtained the nonlinear stability for Couette flow of three dimensional compressible Navier-Stokes equations with Navier boundary conditions on the lower flat boundary and moving (Dirichlet) condition on the upper flat boundary in which the friction coefficient is restricted to be negative.

As pointed out above, it is reasonable to consider the stability and instability problem of the trivial steady state $(v_s=0, p=p_s)$ ($p_s$ is a constant) to equation (\ref{0.1}) under the Navier-slip boundary condition (\ref{0.2}) both for $\alpha\leq 0$ and for $\alpha>0$.

Suppose that ${\mathbf u}$ is a strong solution to the linearized equations perturbed around the steady state $(0,p_s)$ with constant slip length $\alpha$, one deduces that the kinetic energy $\mathcal{E}(t):=\frac{1}{2}\|{\bf u}(t)\|^2_{L^2(\Omega)}$ satisfies the basic energy law
\[\mathcal{E}(t)=\mathcal{E}(0)-\mu\int_0^t\int_{\Omega}|\nabla{\bf u}|^2dxdyds+\alpha\int_0^t\int_{\partial\Omega}|u|^2dxds\]
which indicates that if $\alpha\leq 0$, the total kinetic energy of the system is decay and the stability is not difficult to be proved. For this classical case, we would like to call the boundary condition as {\it dissipative} Navier boundary condition.
However, the situation $\alpha>0$ is much more complicated and the stability or the instability for the steady state $(0,p_s)$ has not been studied yet. In this case, whether the total energy is decay or not depends crucially on the competition between the elasticity energy $\mu\int_0^t\int_{\Omega}|\nabla{\bf u}|^2dxdyds$ and the boundary kinetic energy $\alpha\int_0^t\int_{\partial\Omega}|u|^2dxds$ as we can see in the below sections. At this time, we call the boundary condition as {\it absorptive} Navier boundary condition.

In this paper, for simplicity, we will investigate the following equations in the 2 dimensional slab domain $\Omega=\mathbb{R}\times (0,1)$
\begin{align}\label{1.1}
\begin{cases}
\partial_t{\mathbf v}+{\mathbf v}\cdot\nabla{\mathbf v}+\nabla p-\mu\Delta {\mathbf v}=0,\ \mathrm{in}\ \mathbb{R}\times (0,1), t\geq 0,\\
\mathrm{div}{\mathbf v}=0,\ \mathrm{in}\ \mathbb{R}\times (0,1), t\geq 0;
\end{cases}
\end{align}
with the Navier boundary conditions
\begin{align}
&{\mathbf v}\cdot {\mathbf n}=0,\ \mathrm{on}\ \{y=0,1\},\label{1.2}\\
&[(-p\mathbb{I}+\mu(\nabla {\mathbf v}+\nabla^T{\mathbf v}))\cdot {\mathbf n}]\cdot\tau=k_1{\mathbf v}\cdot\tau,\ \mathrm{on}\ \{y=1\},\label{1.3}\\
&[(-p\mathbb{I}+\mu(\nabla {\mathbf v}+\nabla^T{\mathbf v}))\cdot {\mathbf n}]\cdot\tau=k_0{\mathbf v}\cdot\tau,\ \mathrm{on}\ \{y=0\},\label{1.4}
\end{align}
where the superscript $T$ means matrix transposition, $\mathbb{I}$ is the $2\times2$ identity matrix,  ${\mathbf v(x,y;t)}=(v^1(x,y;t),v^2(x,y;t))$ and $p(x,y;t)$ are the velocity and pressure of the flow respectively, $\mathbf n$ is the outward unit normal vector and $\tau$ is the corresponding tangent vector of the boundary. In our consideration, $\mathbf n=(0,1)$ on $\{y=1\}$ and $\mathbf n=(0,-1)$ on $\{y=0\}$, while $\tau=(1,0)$ on both $\{y=1\}$ and $\{y=0\}$. The viscosity $\mu$ is supposed to be strictly positive and the coefficients $k_0$ and $k_1$ do not have defined sign.

Our main interest here is to study the linear and nonlinear stability and instability of the steady state solution, $(0,p_s)$, to this boundary value problem. Our results show that the stability (or instability) of this equilibrium depends crucially on whether the boundaries dissipate energy and the strengthen of the viscosity and slip length. It is shown, in the case that all the boundaries are dissipative, that nonlinear asymptotic stability holds true. Otherwise, there is a sharp critical viscosity, which distinguishes the nonlinear stability from instability.

Denote the perturbation by
  \[{\mathbf u}={\mathbf v}-{\mathbf 0},\ q=p-p_s.\]

Then $({\mathbf u},q)$ satisfies the perturbed equations
\begin{align}\label{nonl}
\begin{cases}
\partial_t{\mathbf u}+{\mathbf u}\cdot\nabla{\mathbf u}+\nabla q-\mu\Delta {\mathbf u}=0,\\
\mathrm{div}{\mathbf u}=0.
\end{cases}
\end{align}

The boundary conditions, (\ref{1.2})-(\ref{1.4}), can be rewritten as follows
\begin{align}
&u^2(x,0)=u^2(x,1)=0,\ x\in \mathbb{R},\label{bdy1}\\
&\partial_yu^1(x,1)=\frac{k_1}{\mu}u^1(x,1),\ x\in \mathbb{R},\label{bdy2}\\
&\partial_yu^1(x,0)=-\frac{k_0}{\mu}u^1(x,0),\ x\in \mathbb{R}.\label{bdy3}
\end{align}

Linearizing (\ref{nonl}) around the steady state $({\mathbf 0},p_s)$ yields the linearized equations
\begin{align}\label{lin}
\begin{cases}
\partial_t{\mathbf u}+\nabla q-\mu\Delta {\mathbf u}=0,\\
\mathrm{div}{\mathbf u}=0.
\end{cases}
\end{align}

 For convenience, we will use the following notations throughout this paper.
 \[\Omega:=\mathbb{R}\times(0,1),\ L^p:=L^p(\Omega),\ H^k:=W^{2,k}(\Omega), \ \int:=\int_\Omega, \]
 \[H_\sigma^1:=\{\mathbf{u}\in H^1|{\rm div}\mathbf{u}=0, {u}^2=0\ {\rm on}\ \mathbb{R}\times\{0,1\}\}. \]
 $H_0^1(0,1)$ and $H^2(0,1)$ will be written as $H_0^1$ and $H^2$ respectively. In addition, a product space $(X)^2$ of vector functions is still denoted by $X$, for example, the vector function ${\mathbf u}\in (H^1)^2$ is denoted by ${\mathbf u}\in H^1$. $\mathbb{N}$ is the set of nonnegative integers.



First, we study the linear instability of the steady state $({\mathbf 0},p_s)$. To this end, we looks for a growing mode solution to the linearized problem (\ref{bdy1})-(\ref{lin}) in the  form
\begin{align}\label{1.8}
{\mathbf v}(x,y;t)={\mathbf w}(x,y)e^{\lambda t},\ q(x,y;t)={\tilde p}(x,y)e^{\lambda t}
\end{align}
for some $\lambda >0$. Putting this ansatz into (\ref{bdy1})-(\ref{lin}) yields
\begin{align}\label{1.9}
\begin{cases}
\lambda{\mathbf w}+\nabla{\tilde p}-\mu\Delta {\mathbf w}=0,\ \mathrm{in}\ \Omega,\\
\mathrm{div}{\mathbf w}=0,\ \mathrm{in}\ \Omega,
\end{cases}
\end{align}
and the boundary conditions
\begin{align}
&w^2(x,0)=w^2(x,1)=0,\ x\in \mathbb{R};\label{1.10}\\
&\partial_yw^1(x,1)=\frac{k_1}{\mu}w^1(x,1),\ x\in \mathbb{R};\label{1.11}\\
&\partial_yw^1(x,0)=-\frac{k_0}{\mu}w^1(x,0),\ x\in \mathbb{R}.\label{1.12}
\end{align}

We will solve problem (\ref{1.9})-(\ref{1.12}) by the standard normal mode analysis, see \cite{7}. That is, rewrite $\mathbf{w}$ and $\tilde{p}$ in terms of the new unknowns $\phi,\psi,\pi:(0,1)\rightarrow \mathbb{R}$  for each frequency $\xi$ as:
\begin{align}\label{1.13}
w^1(x,y)=-i\phi(y)e^{ix\xi},w^2(x,y)=\psi(y)e^{ix\xi},\tilde{p}(x,y)=\pi(y)e^{ix\xi}.
\end{align}

For each fixed $\xi\neq 0$, this leads to the following system of ODEs
\begin{align}\label{1.14}
\begin{cases}
-\lambda\phi+\xi\pi-\mu\xi^2\phi+\mu\phi^{\prime\prime}=0,\\
\lambda\psi+\pi^\prime+\mu\xi^2\psi-\mu\psi^{\prime\prime}=0,\\
\xi\phi+\psi^\prime=0,
\end{cases}
\end{align}
with boundary conditions
\begin{align}
&\psi(0)=\psi(1)=0,\label{1.15}\\
&\phi^\prime(1)=\frac{k_1}{\mu}\phi(1),\label{1.16}\\
&\phi^\prime(0)=-\frac{k_0}{\mu}\phi(0).\label{1.17}
\end{align}

Eliminating $\pi$ from the second equation of (\ref{1.14}) gives a fourth order ODE for $\psi$
\begin{align}\label{1.18}
-\lambda(\xi^2\psi-\psi^{\prime\prime})=\mu(\psi^{(4)}-2\xi^2\psi^{\prime\prime}+\xi^4\psi),y\in(0,1)
\end{align}
with the boundary conditions
\begin{align}
&\psi(1)=\psi(0)=0,\label{1.19}\\
&\psi^{\prime\prime}(1)=\frac{k_1}{\mu}\psi^\prime(1),\label{1.20}\\
&\psi^{\prime\prime}(0)=-\frac{k_0}{\mu}\psi^\prime(0).\label{1.21}
\end{align}

If, for some frequency $\xi$, there exists a solution to (\ref{1.18})-(\ref{1.21})
with positive $\lambda$, then the above steady state is said to be linearly instable. Since problem (\ref{1.18})-(\ref{1.21}) has a natural variational structure, one may reach such an aim by solving the minimization problem
\begin{align}\label{2.1}
-\lambda=\inf_{H_0^1\cap H^2}\frac{E(\psi)}{J(\psi)},
\end{align}
where
\begin{align}\label{2.2}
E(\psi)=\frac{\mu}{2}\int^1_0\left[(\psi^{\prime\prime})^2+2\xi^2(\psi^\prime)^2+\xi^4\psi^2\right]-\frac{k_1}{2}(\psi^\prime(1))^2
-\frac{k_0}{2}(\psi^\prime(0))^2
\end{align}
and
\begin{align}\label{2.3}
J(\psi)=\frac{1}{2}\int_0^1\left[\xi^2\psi^2+(\psi^\prime)^2\right]
\end{align}
are both well-defined on the space $H_0^1\cap H^2$.

In order to get a positive $\lambda(=\lambda(\xi^2))$ in the variational problem (\ref{2.1}), we observe that if

\[\frac{\mu}{2}\int^1_0(\psi^{\prime\prime})^2-\frac{k_1}{2}(\psi^\prime(1))^2-\frac{k_0}{2}(\psi^\prime(0))^2\]
is negative for small viscosity $\mu$, then $E(\psi)$ is negative for small $\xi$.
This is a key observation which motivates us to define the critical viscosity by
 \begin{align}\label{2.4}
 \mu_c=&\sup_{H_0^1\cap H^2,\psi\neq 0}\frac{k_1(\psi^\prime(1))^2+k_0(\psi^\prime(0))^2}{\int_0^1(\psi^{\prime\prime})^2}.
 \end{align}

It will be shown in next section that explicit values of the critical viscosity are
\begin{equation}\label{2.4111}
 \mu_c=
 \begin{cases}
 0,&\max(k_0,k_0)\leq 0,\\
\displaystyle\frac{k_1+k_0+\sqrt{k_1^2+k_0^2-k_1k_0}}{6},&\max(k_0,k_1)>0.
 \end{cases}
 \end{equation}

Moreover, it will also be shown that this value of $\mu_c$ is a sharp threshold of the stability and instability. Precisely, we have the following main results.

{The first result is on the linear instability.
\begin{thm}\label{thm1}
{\rm (Linear instability)} The steady state $(\mathbf{0},p_s)$ is linearly unstable in $H^k$, for any $k\in\mathbb{N}$, in the sense that there are exponentially growing mode solutions to the linearized perturbed problem (\ref{bdy1})-(\ref{lin}) in $H^k$ if and only if $\mu_c>0$ and $\mu\in(0,\mu_c)$.
\end{thm}

The nonlinear instability is stated in the following theorem.
\begin{thm}\label{thm2}
{\rm (Nonlinear instability)} The steady state $(\mathbf{0},p_s)$ is nonlinearly unstable in $L^2$ norm if and only if $\mu_c>0$ and $\mu\in(0,\mu_c)$. More precisely, we have

(i) Assume $\mu_c>0$ and $\mu\in(0,\mu_c)$. Then there exists a constant $\varepsilon>0$ and a function ${\mathbf u}_0$,  $\|\mathbf{u}_0\|_{H^2}=1$, such that for any $\delta$: $0<\delta<\varepsilon$, there exists a unique global strong solution ($\mathbf{u}^\delta, q$) to the nonlinear perturbed problem (\ref{nonl})-(\ref{bdy3}) with the initial data $\mathbf{u}_0^\delta:=\delta\mathbf{u}_0$, such that ${\mathbf u}^\delta\in C([0,T],H^2)$, $\nabla q\in L^2$ and
\begin{align}\label{2.31111}\|{\mathbf u}^\delta(T^\delta)\|_{L^2}\geq\varepsilon.\end{align}
Here the escape time is $T^\delta:=\frac{1}{\lambda_*}\ln(\varepsilon/\delta)\in(0,T)$, where $\lambda_*$ is defined in (\ref{5.401}).

(ii) Assume $\mu_c\geq 0$ and $\mu\in [\mu_c, +\infty)$. Then there exists a unique global strong solution ($\mathbf{u}, q$) to the nonlinear perturbed problem (\ref{nonl})-(\ref{bdy3}) with ${\mathbf u}\in C([0,T],H^2)$ and $\nabla q\in L^2$, such that
\begin{align}\label{2.32222}\|{\mathbf u}\|_{L^2}\leq \|{\mathbf u}_0\|_{L^2}\end{align}
where ${\mathbf u}_0$ is the initial data of problem (\ref{nonl})-(\ref{bdy3}).

\end{thm}

Finally, we have the following stability theorem for $\mu>\mu_c$, where $\mu_c\geq 0$.
\begin{thm}\label{thm3}
{\rm(Nonlinear asymptotic stability)} The steady states $({\mathbf 0},p_s)$ is nonlinear asymptotically stable globally provided that $\mu>\mu_c\geq 0$, that is, let $(\mathbf{u},q)$ be a solution of nonlinear problem (\ref{nonl})-(\ref{bdy3}) with initial data $\mathbf{u}_0$, then the followings are true.

(i) For general initial data $\mathbf{u}_0$, there exist some constants $\alpha>0, C>0$, such that
 \begin{align}\label{6.01}
 \sup_{0<t<+\infty}\|\mathbf{u}(t)\|_{H^2}\leq C\|\mathbf{u}_0\|_{H^2},
 \lim_{t\rightarrow+\infty}\|\mathbf{u}(t)\|_{H^2}=0,\mathrm{~and~}
 \|\mathbf{u}(t)\|_{H^1}\leq Ce^{-\alpha t},
 \end{align}
 where $C=C(\mu,k_0,k_1,\|\mathbf{u}_0\|_{H^2})$ is increasing with respect to $\|\mathbf{u}_0\|_{H^2}.$

 (ii) Moreover, if the initial data is small, then we have $H^2$-norm decay estimates. That is, there exists a positive constant $\beta$ such that
\begin{align}\label{6.2}
\|{\mathbf u}(t)\|_{H^2}\leq Ce^{-\beta t},
\end{align}
 provided that the initial data ${\mathbf u}_0$ satisfies $\|{\mathbf u}_0\|_{H^2}\leq\bar{\sigma}$, for some constant $\bar{\sigma}>0$.
\end{thm}}

{The rest of this paper is arranged as follows. First, we analyse in detail the problem (\ref{2.4}) to determine the exact value of the critical viscosity. In Section 3 and Section 4, we prove the instability part of Theorem \ref{thm1} and Theorem \ref{thm2} respectively. The stability part of Theorem \ref{thm1}, Theorem \ref{thm2} and Theorem \ref{thm3} will be proved in Section 5.}

\section{The critical viscosity}
In order to obtain the value of the critical viscosity, we consider the equivalent variational problem of (\ref{2.4}) as
\begin{equation}\label{2.8000}
\sup_{\psi\in\mathcal{Y}}Z(\psi),
\end{equation}
where
\begin{align}
&Z(\psi)=\frac{k_1}{2}(\psi^\prime(1))^2+\frac{k_0}{2}(\psi^\prime(0))^2,\label{2.5}\\
&\mathcal{Y}=\left\{\psi\in H_0^1\cap H^2\Bigg|\frac{1}{2}\int_0^1(\psi^{\prime\prime})^2=1\right\}.\label{2.6}
\end{align}

In what follows, we shall use the fact that, for any $f\in H_0^1\cap H^2$, there holds
\begin{align}\label{2.8}
\|f^\prime\|^2_{L^2}\leq\|f^{\prime\prime}\|^2_{L^2},
\end{align}
where one can use Poincar\'e inequality and the fact $f\in H^1_0$ to prove (\ref{2.8}).

Thus, for any $\psi\in\mathcal{Y}$,
\begin{align}\label{2.9}
|Z(\psi)|=&\frac{1}{2}\left|\int_0^1\left[((k_1+k_0)y-k_0)(\psi^\prime)^2\right]^\prime_y dy\right|\nonumber\\
=&\left|\frac{(k_1+k_0)}{2}\int_0^1(\psi^\prime)^2dy+\int_0^1[(k_1+k_0)y-k_0]\psi^\prime\psi^{\prime\prime}dy\right|\nonumber\\
\leq&\frac{1}{2}\int_0^1\left[ C_1(\psi^\prime)^2+C_2(\psi^{\prime\prime})^2\right]dy\nonumber\\
\leq& \frac{1}{2}\int_0^1C_3(\psi^{\prime\prime})^2dy=C_3,
\end{align}
for some positive constants $C_1,C_2$ and $C_3$, depending only on $k_1$ and $k_0$. This shows that $\sup_{\psi\in\mathcal{Y}}Z(\psi)$ exists and is finite.

 Set $\mu_c:=\sup_{\psi\in\mathcal{Y}}Z(\psi)$. The exact values of $\mu_c$ will be given in different cases in the following two propositions.
 \begin{prop}\label{p2.1}
 Let $\mu_c\in\mathbb{R}$ be defined in (\ref{2.4}). Then $\mu_c=0$ for $\max(k_0,k_1)\leq 0$ and $\mu_c>0$ while $\max(k_0,k_1)>0$.
 \end{prop}
 \proof
 If both $k_0$ and $k_1$ are non-positive, then, clearly,
 \begin{align}\label{2.91}
 \mu_c=\sup_{\psi\in\mathcal{Y}}\left[\frac{k_1}{2}(\psi(1))^2+\frac{k_0}{2}(\psi(0))^2\right]\leq 0.
 \end{align}
 On the other hand, for a suitable choice of $\alpha$,
 \begin{equation}\label{2.92}
 \psi_1(x):=
 \begin{cases}
 0,& x\in[0,\frac{1}{4}];\\
 \alpha\exp\{\frac{1}{(x-\frac{1}{2})^2-\frac{1}{16}}\},& x\in(\frac{1}{4},\frac{3}{4});\\
 0,& x\in[\frac{3}{4},1]
 \end{cases}
 \end{equation}
 belongs to $\mathcal{Y}$. Moreover, $\psi_1\in C_0^\infty([0,1])$ and $\psi^\prime_1(1)=\psi_1^\prime(0)=0$, which implies $Z(\psi_1)=0$. This,
 together with (\ref{2.91}), implies that $\mu_c=0$.

 In the other case, without loss of generality, we suppose that $k_0>0$ and define
 \begin{equation}\label{2.93}
 \psi_2(x):=
 \begin{cases}
 -\frac{2}{3}x(x-\frac{3\sqrt{2}-3}{2}),&x\in[0,\frac{3\sqrt{2}}{8}];\\
 \frac{2}{3}(x-\frac{3}{4})^2,&x\in(\frac{3\sqrt{2}}{8},\frac{3}{4});\\
 0,&x\in[\frac{3}{4},1].
 \end{cases}
 \end{equation}
 Then, one can verify that $\psi_2\in\mathcal{Y}$ and $Z(\psi_2)>0$, which means $\mu_c>0$ in this case.
 \endproof

 To find the exact value of $\mu_c$ in the case of $\mu_c>0$, we need the following Proposition.
 \begin{prop}\label{p2.2}
 Let $\mu_c$ be defined as in (\ref{2.4}) and suppose that $\max\{k_0,k_1\}>0$. Then
 \begin{equation*}
 \mu_c=\frac{k_1+k_0+\sqrt{k_1^2+k_0^2-k_1k_0}}{6}.
 \end{equation*}
 \end{prop}
 \proof
Let $\{\psi_n\}_{n=1}^\infty\in\mathcal{Y}$ be a maximizing sequence. It follows from (\ref{2.8}) that
\[\|\psi_n\|_{H_0^1\cap H^2}=\|\psi^\prime_n\|_{L^2}+\|\psi^{\prime\prime}_n\|_{L^2}\leq C_4 \]
for some constant $C_4>0$. Therefore $\{\psi_n\}$ is bounded in $H_0^1\cap H^2$. Hence, up to a subsequence if necessary, $\psi_n\rightharpoonup \psi$ weakly in $H^2$ and $\psi_n\rightarrow \psi$ strongly in $H^1_0$. This implies that
\begin{align}\label{2.10}
Z(\psi)&=\frac{k_1+k_0}{2}\int_0^1(\psi^\prime)^2+\int_0^1[(k_1+k_0)y-k_0]\psi^\prime\psi^{\prime\prime}\nonumber\\
&=\frac{k_1+k_0}{2}\lim_{n\rightarrow\infty}\int_0^1(\psi_n^\prime)^2+\lim_{n\rightarrow\infty}\int_0^1[(k_1+k_0)y-k_0]\psi_n^\prime\psi_n^{\prime\prime}\nonumber\\
&=\lim_{n\rightarrow\infty}Z(\psi_n)=\mu_c,
\end{align}
and
\begin{align}\label{2.11}
\frac{1}{2}\|\psi^{\prime\prime}\|^2_{L^2}\leq\lim_{n\rightarrow\infty}\frac{1}{2}\|\psi_n^{\prime\prime}\|^2_{L^2}=1.
\end{align}

Now we claim that $\|\psi^{\prime\prime}\|^2_{L^2}=2$, i.e. $\psi\in\mathcal{Y}$.

Otherwise, one may assume that $\|\psi^{\prime\prime}\|^2_{L^2}:=2r^2<2$, for some constant $0\leq r<1$. Notice that if $r=0$, then $\psi=0$, which implies that $\mu_c=0$. Thus, $0<r<1$ and $\tilde{\psi}=\psi/r\in\mathcal{Y}$. The definition of $\mu_c$ and (\ref{2.10}) lead to
$$
\mu_c\geq Z(\tilde{\psi})=Z(\psi)/r^2=\mu_c/r^2>\mu_c,
$$
which is a contradiction. Thus, $\psi\in\mathcal{Y}$ is a maximizer of the variational problem (\ref{2.8000}).

In what follows, we will find the exact expression of the unique maximizer $\psi$ and then obtain the exact value of $\mu_c$.

For any $\psi_0\in H_0^1\cap H^2$, $s,r\in\mathbb{R}$, define
\begin{align}\label{2.12}
I(s,r)=\frac{1}{2}\int_0^1(\psi^{\prime\prime}+s\psi_0^{\prime\prime}+r\psi^{\prime\prime})^2.
\end{align}

Notice that $I(s,r)$ is smooth and that
\begin{align}\label{2.13}
I(0,0)=1,
\end{align}
\begin{align}
&\partial_sI(s,r)\big|_{(0,0)}=\int_0^1(\psi^{\prime\prime}+s\psi_0^{\prime\prime}+r\psi^{\prime\prime})\psi_0^{\prime\prime}\Big|_{(0,0)}
=\int_0^1\psi^{\prime\prime}\psi_0^{\prime\prime},\label{2.14}\\
&\partial_rI(s,r)\big|_{(0,0)}=\int_0^1(\psi^{\prime\prime}+s\psi_0^{\prime\prime}+r\psi^{\prime\prime})\psi^{\prime\prime}\Big|_{(0,0)}
=\int_0^1(\psi^{\prime\prime})^2=2\neq 0.\label{2.15}
\end{align}

By implicit function theorem, there exists a smooth function $r=r(s)$ defined near $s=0$ such that $r(0)=0,\ I(s,r(s))\equiv 1.$
It follows from this and the fact that $\psi$ is a maximizer that
\begin{align}\label{2.16}
0&=\frac{d}{ds}\Bigg|_{s=0}Z(\psi+s\psi_0+r(s)\psi)\nonumber\\
&=k_1\psi^\prime(1)\psi_0^\prime(1)+k_0\psi^\prime(0)\psi_0^\prime(0)+r^\prime(0)(k_1(\psi^\prime(1))^2+k_0(\psi^\prime(0))^2)
\end{align}
for any test function $\psi_0\in H^1_0\cap H^2$.

Differentiating the equation $I(s,r(s))=1$ yields that
\begin{align}\label{2.17}
I^\prime(s,r(s))=\partial_sI(s,r(s))+\partial_rI(s,r(s))r^\prime(s)\equiv 0,
\end{align}
which implies that
\begin{align}\label{2.18}
r^\prime(0)=-\frac{1}{2}\int\psi^{\prime\prime}\psi_0^{\prime\prime}.
\end{align}

It follows from (\ref{2.18}) and (\ref{2.16}) that
\begin{align}\label{2.19}
\mu_c\int_0^1\psi^{\prime\prime}\psi_0^{\prime\prime}=k_1\psi^\prime(1)\psi_0^\prime(1)+k_0\psi^\prime(0)\psi_0^\prime(0).
\end{align}

First, choosing $\psi_0$ to be compactly supported in $(0,1)$ in (\ref{2.19}) shows that
 \begin{align}\label{2.20}
 \mu_c\psi^{(4)}=0,x\in(0,1),
 \end{align}
 in a weak sense. Standard bootstrapping arguments show that the solution $\psi$ is smooth. Then (\ref{2.20}) and (\ref{2.19}) implies that
\begin{align*}
&\mu_c\psi^{\prime\prime}(1)=k_1\psi^\prime(1),\\
&\mu_c\psi^{\prime\prime}(0)=-k_0\psi^\prime(0).
\end{align*}

Therefore, the maximizers $\psi$ must solve the following problem
\begin{align}\label{2.21}
\begin{cases}
\mu_c\psi^{(4)}=0, y\in(0,1);\\
\mu_c\psi^{\prime\prime}(1)=k_1\psi^\prime(1),\\
\mu_c\psi^{\prime\prime}(0)=-k_0\psi^\prime(0),\\
\psi(1)=\psi(0)=0.
\end{cases}
\end{align}

In view of $(\ref{2.21})_1$ and $(\ref{2.21})_4$, we set that
\begin{align}\label{2.22}
\psi(x)=x(x-1)(Ax+B),
\end{align}
where $A,B$ are coefficients to be determined. It is easy to see that
\begin{align*}
&\psi'(x)=Ax(x-1)+(x-1)(Ax+B)+x(Ax+B),\\
&\psi''(x)=6Ax-2A+2B.
\end{align*}
Then, we get
\begin{align}
&\psi'(1)=A+B,~~\psi'(0)=-B;\label{2.23}\\
&\psi''(1)=4A+2B,~~\psi''(0)=-2A+2B.\label{2.24}
\end{align}
Substituting \eqref{2.23} and \eqref{2.24} into $\eqref{2.21}_2$ and $\eqref{2.21}_3$ gives
\begin{align}\label{2.25}
\begin{cases}
\mu_c(4A+2B)=k_1(A+B),\\
\mu_c(-2A+2B)=k_0B,
\end{cases}
\end{align}
which can be rewritten as
\begin{align}\label{2.26}
\begin{cases}
(4\mu_c-k_1)A+(2\mu_c-k_1)B=0,\\
-2\mu_cA+(2\mu_c-k_0)B=0.
\end{cases}
\end{align}
By the theory of linear algebra, the above linear equations admit a non-zero solution only if the determinant of the coefficient matrix equals to 0, that is,
\begin{align}\label{2.27}
\left|
\begin{array}{cc}
  4\mu_c-k_1,&2\mu_c-k_1\\
  -2\mu_c,&2\mu_c-k_0
\end{array}
\right|=0.
\end{align}

This is equivalent to
\[12\mu_c^2-4\mu_c(k_0+k_1)+k_0k_1=0,\]
which yields
\begin{align}\label{2.28}
\mu_c=\frac{k_0+k_1\pm\sqrt{k_0^2+k_1^2-k_0k_1}}{6}.
\end{align}
Since $\mu_c$ is the maximum of the corresponding variational problem, we should take
\begin{align}\label{2.29}
\mu_c=\frac{k_0+k_1+\sqrt{k_0^2+k_1^2-k_0k_1}}{6}.
\end{align}

In addition, substituting the expression of $\mu_c$ into \eqref{2.26} gives the relationship of $A$ and $B$, and the unique maximizer is also given, with the coefficients uniquely determined by using $\int_0^1(\psi^{\prime\prime})^2=2.$
\endproof

\begin{rem}\label{r2.3}
As was pointed out in \cite{TTN} with a counterexample, $\mu_c=\frac{k}{6}$ for the case that $k_0=k_1=k>0$ in the previous version(arXiv:1608.03019) of this paper should be revised. The origin of this mistake in the previous version is that, compared with the undetermined maximizer given by \eqref{2.22} in this version, we set that in the form $\psi(x)=ax(x-1)(x-b)$, which restricts that $a\neq 0.$ In this version, we have corrected this error. In addition, if we consider a more general domain $\mathbb{R}\times(\alpha,\beta)$ with $\alpha<\beta$, to investigate the effect of the width $\beta-\alpha$ on the value of $\mu_c$, we find that
\[\mu_c=\frac{(\beta-\alpha)(k_0+k_1)+\sqrt{(k_0^2+k_1^2-k_0k_1)(\beta-\alpha)^2}}{6},~~\text{for}~\max(k_0,k_1)>0.\]
When $\alpha=-1$ and $\beta=1$, this formula is consistent with the corresponding result in \cite{TTN}, where the author obtained that via Lagrangian multiplier method.
\end{rem}

\section{The linear instability}
\subsection{Analysis for the variational problem (\ref{2.1})}
In order to prove the instability part of Theorem \ref{thm1} in this section, we will discuss the corresponding variational problem (\ref{2.1})-(\ref{2.3}) with a fixed $\xi$ by variational methods.

To rewrite this variational problem in an equivalent form, we define
\begin{align}\label{3.1}
\mathcal{A}=\{\psi\in H_0^1((0,1))\cap H^2((0,1))|J(\psi)=1\}
\end{align}
where $J(\psi)$ is defined by (\ref{2.3}).

The main task is to show that the minimum of $E(\psi)$ over $\mathcal{A}$ can be achieved and the minimizer solves the Euler-Lagrange equation, which is equivalent to (\ref{1.18}) together with the corresponding boundary conditions (\ref{1.19})-(\ref{1.21}).
First, the existence of the minimizer is shown below.
\begin{prop}\label{p3.1}
$E(\psi)$ achieves its minimum on $\mathcal{A}$.
\end{prop}
\proof
Using the constraint on $J(\psi)$ and Cauchy inequality, we get
\begin{align}\label{3.2}
E(\psi)&=\frac{\mu}{2}\int_0^1\left[(\psi^{\prime\prime})^2+2\xi^2(\psi^\prime)^2+\xi^4\psi^2\right]
-\frac{1}{2}\int_0^1\left[\left((k_1+k_0)y-k_0\right)(\psi^\prime)^2\right]^\prime_ydy\nonumber\\
&=\frac{\mu}{2}\int_0^1(\psi^{\prime\prime})^2+\frac{\mu\xi^2}{2}\int_0^1\left[2(\psi^\prime)^2+\xi^2\psi^2\right]\nonumber\\
&\quad-\frac{1}{2}\int_0^1\left[(k_1+k_0)(\psi^\prime)^2+2\psi^\prime\psi^{\prime\prime}\left((k_1+k_0)y-k_0\right)\right]\nonumber\\
&\geq\mu\xi^2-\frac{1}{2}\int^1_0\left[(k_1+k_0)+\mu^{-1}\left((k_1+k_0)y-k_0\right)^2\right](\psi^\prime)^2\nonumber\\
&\geq\mu\xi^2-\frac{C_0}{2}\int_0^1(\psi^\prime)^2\geq\mu\xi^2-C_0
\end{align}
for any fixed $\xi\in \mathbb{R}$, where
\begin{align}\label{3.3}
C_0=\max_{0\leq y\leq 1}\left[|k_1+k_0|+\mu^{-1}\left((k_1+k_0)y-k_0\right)^2\right].
\end{align}
This means that $E$ is bounded from below over $\mathcal{A}$, and thus $\inf_\mathcal{A}E(\psi)$ is well defined and finite.

Denote $-\lambda :=\inf_\mathcal{A}E(\psi)$, and let $\{\psi_n\}^\infty_{n=1}\in\mathcal{A}$ be a minimizing sequence. Without lose of generality, one may assume that $E(\psi_n)\leq -\lambda +1$. Then the constraint on $J(\psi_n)$ and the Poincar\'e inequality imply that ${\psi_n}$ is uniformly bounded in $H^1$, which is independent of $\xi^2$. In addition, by the definition of $E$ and the Cauchy inequality, one has
\begin{align}\label{3.5}
\mu\int_0^1(\psi_n^{\prime\prime})^2&\leq 2E(\psi_n)+\int_0^1\left[\left((k_1+k_0)y-k_0\right)(\psi_n^\prime)^2\right]^\prime\nonumber\\
&=2E(\psi_n)+\int_0^1\left[(k_1+k_0)(\psi_n^\prime)^2+2\psi_n^\prime\psi_n^{\prime\prime}\left((k_1+k_0)y-k_0\right)\right]\nonumber\\
&\leq 2E(\psi_n)+\frac{\mu}{2}\int_0^1(\psi_n^{\prime\prime})^2+2C_0\int^1_0(\psi_n^\prime)^2,
\end{align}
which implies that
\begin{align}\label{3.6}
\int_0^1(\psi_n^{\prime\prime})^2\leq 4\mu^{-1}(E(\psi_n)+C_0).
\end{align}

It follows that the sequence $\{\psi_n\}$ is bounded in $H_0^1\cap H^2$, and thus, up to a subsequence if necessary, $\psi_n\rightharpoonup\psi$ weakly in $H^2$ and $\psi_n\rightarrow\psi$ strongly in $H_0^1$.

Rewrite $E$ as
\begin{align}\label{3.7}
E(\psi)=\frac{\mu}{2}\int_0^1\left[(\psi^{\prime\prime})^2+\left(\xi^2-\mu^{-1}(k_0+k_1)\right)(\psi^\prime)^2+\xi^4\psi^2\right]
-\int_0^1\psi^\prime\psi^{\prime\prime}\left((k_1+k_0)y-k_0\right).
\end{align}

It follows from the weak lower semi-continuity and weak convergence in $H^2$ and strong convergence in $H^1_0$ that
\begin{align}\label{3.8}
E(\psi)&\leq\frac{\mu}{2}\liminf_{n\rightarrow\infty}\int_0^1(\psi_n^{\prime\prime})^2+\left(\mu\xi^2-\frac{(k_0+k_1)}{2}\right)
\lim_{n\rightarrow\infty}\int_0^1(\psi_n^\prime)^2\nonumber\\
&\qquad+\frac{\mu\xi^4}{2}\lim_{n\rightarrow\infty}\int_0^1\psi_n^2
-\lim_{n\rightarrow\infty}\int_0^1\psi_n^\prime\psi_n^{\prime\prime}\left((k_1+k_0)y-k_0\right)\nonumber\\
&=\liminf_{n\rightarrow\infty}E(\psi_n)=\inf_\mathcal{A}E(\psi).
\end{align}

Finally, the claim $J(\psi)=1$ follows from the strong convergence in $H^1_0$.
\endproof
\begin{rem}\label{r3.1}
{\it Since the aim here is to look for growing mode solutions to the linearized equation (\ref{lin}) with boundary conditions (\ref{bdy1})-(\ref{bdy3}), one should restrict the parameter $\xi$ to stay in a specific range to guarantee that $\inf_\mathcal{A}E(\psi)=-\lambda<0$. It requires, in view of (\ref{3.2}), at least that
\[C_0>0\ \mathrm{and}\ \xi^2\leq \mu^{-1} C_0,\]
so that it is possible to have a negative minimum for $E(\psi)$ over $\mathcal{A}$. This will be achieved later by finding a critical frequency.}
\end{rem}

Next we will show that the minimizer constructed above satisfies an Euler-Lagrangian equation equivalent to (\ref{1.18}).
\begin{prop}\label{p3.2}
Let $\psi\in\mathcal{A}$ be the minimizer of $E$ constructed in {\it Proposition \ref{p3.1}}, and denote $-\lambda:=E(\psi)$. Then $\psi$ is smooth and satisfies
\begin{align}\label{3.9}
-\lambda(\xi^2\psi-\psi^{\prime\prime})=\mu(\psi^{(4)}-2\xi^2\psi^{\prime\prime}+\xi^4\psi),
\end{align}
along with the boundary conditions
\begin{align}
&\psi(1)=\psi(0)=0,\label{3.10}\\
&\psi^{\prime\prime}(1)=\frac{k_1}{\mu}\psi^\prime(1),\label{3.11}\\
&\psi^{\prime\prime}(0)=-\frac{k_0}{\mu}\psi^\prime(0).\label{3.12}
\end{align}
As a consequence, there exists a solution $(\phi,\psi,\pi)$ to the problem (\ref{1.14})-(\ref{1.17}).
\end{prop}
\proof
For any $\psi_0\in H_0^1\cap H^2,\ t,r\in \mathbb{R}$, let $\psi\in \mathcal{A}$ be a minimizer and define
 \[j(t,r):=J(\psi+t\psi_0+r\psi).\] Then $j(t,r)$ is smooth and $j(0,0)=1,$. Notice that
\begin{align}\label{3.13}
\partial_tj(0,0)=\int_0^1(\xi^2\psi\psi_0+\psi^\prime\psi_0^\prime),\ \mathrm{and}\ \partial_rj(0,0)=\int_0^1(\xi^2\psi^2+(\psi^\prime)^2)=2\neq 0.
\end{align}

 Then, by implicit function theorem, there exists a smooth function $r=r(t)$ defined near $0$ such that $r(0)=0$ and $j(t,r(t))=1$.

 Since $\psi$ is a minimizer, it is clear that
 \begin{align}\label{3.14}
 0&=\frac{\mathrm{d}}{\mathrm{d}t}\bigg|_{t=0}E(\psi+t\psi_0+r(t)\psi)\nonumber\\
 &=\mu\int_0^1\left(\psi^{\prime\prime}\psi_0^{\prime\prime}+2\xi^2\psi^\prime\psi_0^\prime+\xi^4\psi\psi_0\right)
 -k_1\psi^\prime(1)\psi_0^\prime(1)-k_0\psi^\prime(0)\psi_0^\prime(0)+2r^\prime(0)E(\psi).
 \end{align}

 Now differentiating the equation $j(t,r(t))=1$ gives
 \begin{align}\label{3.16}
 r^\prime(0)=-\frac{1}{2}\int_0^1(\xi^2\psi\psi_0+\psi^\prime\psi_0^\prime).
 \end{align}

 Substituting (\ref{3.16}) into (\ref{3.14}) yields
 \begin{align}\label{3.17}
 &\mu\int_0^1\left(\psi^{\prime\prime}\psi_0^{\prime\prime}+2\xi^2\psi^\prime\psi_0^\prime+\xi^4\psi\psi_0\right)
 +\lambda\int_0^1(\xi^2\psi\psi_0+\psi^\prime\psi_0^\prime)\nonumber\\
 =&k_1\psi^\prime(1)\psi_0^\prime(1)+k_0\psi^\prime(0)\psi_0^\prime(0).
 \end{align}

 Choosing $\psi_0$ to be compactly supported in $(0,1)$ shows that $\psi_0^\prime(1)=\psi_0^\prime(0)=0$. Substituting this into $(\ref{3.17})$ yields that $\psi\in H^2$ solves (\ref{3.9}) in a weak sense.
 Standard bootstrap arguments then show that the solution is smooth. Next, using equation (\ref{3.9}) and integrating by part lead to
 \begin{align}\label{3.18}
 \mu\int_0^1(\psi^{\prime\prime}\psi_0^\prime)^\prime=k_1\psi^\prime(1)\psi_0^\prime(1)+k_0\psi^\prime(0)\psi_0^\prime(0),
 \end{align}
 which is equivalent to
 \begin{align}\label{3.19}
 \left(\mu\psi^{\prime\prime}(1)-k_1\psi^\prime(1)\right)\psi_0^\prime(1)=\left(\mu\psi^{\prime\prime}(0)+k_0\psi^\prime(0)\right)\psi_0^\prime(0).
 \end{align}

 Since $\psi_0$ is arbitrarily chosen, it follows that
 \begin{align}\label{3.20}
 \psi^{\prime\prime}(1)=\frac{k_1}{\mu}\psi^\prime(1)~~\mathrm{and}~~\psi^{\prime\prime}(0)=-\frac{k_0}{\mu}\psi^\prime(0).
 \end{align}

The Proposition follows. \endproof
\begin{rem}\label{r3.2}
{\it It should be noted that for fixed $\xi^2$, the existence of the solutions $(\phi,\psi,\pi)$ and the corresponding eigenvalue $\lambda$ of problem (\ref{1.14})-(\ref{1.17}) are independent of the values of $k_1,k_0$ and $\mu>0$. That is, for any fixed $\xi^2\in[0,+\infty)$, the functions $(\phi(\xi^2,y),\psi(\xi^2,y),\pi(\xi^2,y))$ and eigenvalue $\lambda(\xi^2)$ are well-defined.}
\end{rem}

To study the sign of $-\lambda(\xi^2)$, which determines the linearized stability of the steady states, we will study the relations among $k_1,k_0,\mu>0$ and $\xi^2$ in details later.

\subsection{Proof of the instability part of Theorem \ref{thm1}}
It follows from the definition of $\mu_c$ that when $\mu_c>0$ and $\mu\in(0,\mu_c)$, there exists $\tilde{\psi}\in H_0^1\cap H^2$, such that
\[\mu\int_0^1(\tilde{\psi}^{\prime\prime})^2-k_1(\tilde{\psi}^\prime(1))^2-k_0(\tilde{\psi}^\prime(0))^2<0,\]

In order to prove the existence of growing mode solutions in this case, it suffices to prove that there is an eigenvalue $\lambda>0$. To do this, since $E(\psi)$ is bounded from below over $\mathcal{A}$, one needs to prove that there exists a function $\tilde{\psi}$ belonging to $\mathcal{A}$ such that $E(\tilde{\psi})<0$.

{\it Step 1.} In this step, we intend to show that there exists  $\tilde{\psi}\in H_0^1\cap H^2$ such that $E(\tilde{\psi})<0$ for some frequency $\xi$, i.e.,
\begin{align}\label{4.3}
\xi^2<\frac{k_1(\tilde{\psi}^\prime(1))^2+k_0(\tilde{\psi}^\prime(0))^2-\mu\int_0^1(\tilde{\psi}^{\prime\prime})^2}
{\mu\int_0^1\left(2(\tilde{\psi}^\prime)^2+\xi^2\tilde{\psi}^2\right)}.
\end{align}

The appearance of $\xi^2$ on the both sides of (\ref{4.3}) makes it difficult to use variational techniques to express the critical value of $\xi^2$. In order to circumvent this difficulty, one can replace the $\xi^2$ on the right-hand side of (\ref{4.3}) with an arbitrary parameter $s^2\geq0$. Precisely, we introduce a family of modified variational problems given by
\begin{align}\label{4.4}
\mathcal{N}^*(s^2)=\sup_{H_0^1\cap H^2,\psi\neq0}\mathcal{N}(\psi,s^2),
\end{align}
where
\begin{align}\label{4.5}
\mathcal{N}(\psi,s^2):=\frac{k_1(\psi^\prime(1))^2+k_0(\psi^\prime(0))^2-\mu\int_0^1(\psi^{\prime\prime})^2}
{\mu\int_0^1\left(2(\psi^\prime)^2+s^2\psi^2\right)},s^2\in[0,+\infty).
\end{align}

Similar to the proof of Proposition \ref{p3.1}, one can prove that $\mathcal{N}^*(s^2)$ is well-defined and the maximizer is achievable for any fixed $s^2\in [0,+\infty)$. Moreover, if $\mu_c>0$ and $\mu\in(0,\mu_c)$, then $\mathcal{N}^*(s^2)>0$ for any $s^2\in[0,+\infty)$.

To establish the continuity, boundedness and monotonicity for the function $\mathcal{N}^*(s^2)$, one sets, for convenience, that
\begin{align}
&\mathcal{N}_1(\psi)=k_1(\psi^\prime(1))^2+k_0(\psi^\prime(0))^2-\mu\int_0^1(\psi^{\prime\prime})^2,\label{4.6}\\
&\mathcal{N}_2(\psi,s^2)=\mu\int_0^1\left(2(\psi^\prime)^2+s^2\psi^2\right).\label{4.7}
\end{align}
\begin{prop}\label{p4.1}
Let $\mathcal{N}^*(s^2):[0,+\infty)\rightarrow\mathbb{R}^+$ be defined by (\ref{4.4})-(\ref{4.5}). Then it holds that

(i) $\mathcal{N}^*(s^2)$ is strictly decreasing,

(ii) $\mathcal{N}^*(s^2)\in C^{0,1}([0,+\infty))$, in particular, $\mathcal{N}^*(s^2)\in C^0([0,+\infty))$.

\end{prop}
\proof
For any $s_1^2,s_2^2\in[0,+\infty)$, define
\[\mathcal{N}^*(s_1^2)=\mathcal{N}(\psi_{s^2_1},s_1^2),\mathcal{N}^*(s_2^2)=\mathcal{N}(\psi_{s^2_2},s_2^2).\]

Then, for $s_1^2<s_2^2$, the definition of supremum and the monotonicity of $\mathcal{N}$ with respect to $s^2$ give
\begin{align}\label{4.11}
\mathcal{N}^*(s_1^2)=\mathcal{N}(\psi_{s^2_1},s_1^2)\geq\mathcal{N}(\psi_{s^2_2},s_1^2),
>\mathcal{N}(\psi_{s^2_2},s_2^2)=\mathcal{N}^*(s_2^2)
\end{align}
which means that $\mathcal{N}^*$ is strictly decreasing with respect to $s^2$, this proves {\it (i)}.

Next, for any $s_1^2,s_2^2\in[0,+\infty)$, by the definition of $\mathcal{N}^*(s^2)$, we have
\begin{align}\label{4.8}
\mathcal{N}^*(s_1^2)&=\frac{\mathcal{N}_1(\psi_{s_1^2})}{\mathcal{N}_2(\psi_{s_1^2},s_1^2)}
=\frac{\mathcal{N}_1(\psi_{s_1^2})}{\mathcal{N}_2(\psi_{s_1^2},s_2^2)}
+\frac{\mathcal{N}_1(\psi_{s_1^2})}{\mathcal{N}_2(\psi_{s_1^2},s_1^2)}
-\frac{\mathcal{N}_1(\psi_{s_1^2})}{\mathcal{N}_2(\psi_{s_1^2},s_2^2)}\nonumber\\
&\leq\mathcal{N}^*(s_2^2)+\frac{\mathcal{N}_1(\psi_{s_1^2})(s_2^2-s_1^2)\int_0^1\psi_{s_1^2}^2}
{\mathcal{N}_2(\psi_{s_1^2},s_1^2)\mathcal{N}_2(\psi_{s_1^2},s_2^2)}
\end{align}

In view of the fact that
\[0<\mathcal{N}_1(\psi_{s_1^2})=C_0\int_0^1(\psi^\prime_{s^2_1})^2\]
and applying Poincar\'e inequality, one gets
\begin{align}\label{4.9}
\mathcal{N}^*(s_1^2)\leq\mathcal{N}^*(s_2^2)+\frac{C_0|s_2^2-s_1^2|}{4\mu^2}=\mathcal{N}^*(s_2^2)+K|s_2^2-s_1^2|,
\end{align}
where $K:=\frac{C_0}{4\mu^2}$. This implies that $\mathcal{N}^*(s^2)\in C^{0,1}([0,+\infty))$.

In addition, for any $s^2\in[0,+\infty)$, it follows from a similar  proof as for (\ref{4.9}) that
\begin{align}\label{4.10}
0<\mathcal{N}^*(s^2)=\sup_{H_0^1\cap H^2,\psi\neq0}\frac{\mathcal{N}_1(\psi)}{\mathcal{N}_2(\psi,s^2)}\leq\frac{C_0}{2\mu}.
\end{align}

This verifies {\it (ii)} and thus the Proposition follows.
\endproof

Now define a function $\Phi:(0,+\infty)\rightarrow(0,+\infty)$ by \[\Phi(s^2)=\frac{s^2}{\mathcal{N}^*(s^2)}.\]

It follows from the properties of $\mathcal{N}^*(s^2)$ that $\Phi$ is continuous and strictly increasing with respect to $s^2$. Since
$\lim_{s^2\rightarrow 0^+}\mathcal{N}^*(s^2)=\mathcal{N}^*(0)>0$ and (\ref{4.10}), thus
\[\lim_{s^2\rightarrow 0^+}\Phi(s^2)=0,\ \mathrm{and}\ \lim_{s^2\rightarrow +\infty}\Phi(s^2)=+\infty.\]

Then by the mean value theorem, there exists $s_0^2\in(0,+\infty)$ such that $\Phi(s_0^2)=1$, i.e., $s_0^2=\mathcal{N}^*(s_0^2)$. Taking $\xi_c^2=s^2_0$ yields that
\begin{align}\label{4.12}
\xi_c^2=\sup_{H_0^1\cap H^2,\psi\neq0}\frac{k_1(\psi^\prime(1))^2+k_0(\psi^\prime(0))^2-\mu\int_0^1(\psi^{\prime\prime})^2}
{\mu\int_0^1\left(2(\psi^\prime)^2+\xi_c^2\psi^2\right)},
\end{align}
which implies that for any $\xi^2\in[0,\xi^2_c)$, it holds that
\begin{align}\label{4.13}
\xi^2<\xi_c^2=\sup_{H_0^1\cap H^2,\psi\neq0}\frac{\mathcal{N}_1(\psi)}{\mathcal{N}_2(\psi,\xi_c^2)}<\sup_{H_0^1\cap H^2,\psi\neq0}\frac{\mathcal{N}_1(\psi)}{\mathcal{N}_2(\psi,\xi^2)}
\end{align}

By the definition of supremum, for any $\xi^2\in [0,\xi^2_c)$, there exists $\tilde{\psi}\in H_0^1\cap H^2$, so that
\[\xi^2<\frac{\mathcal{N}_1(\tilde{\psi})}{\mathcal{N}_2(\tilde{\psi},\xi^2)}=\frac{k_1(\tilde{\psi}^\prime(1))^2+k_0(\tilde{\psi}^\prime(0))^2-\mu\int_0^1(\tilde{\psi}^{\prime\prime})^2}
{\mu\int_0^1\left(2(\tilde{\psi}^\prime)^2+\xi^2\tilde{\psi}^2\right)}.\]
Thus, (\ref{4.3}) is proved.

In order to emphasize the dependence on $\xi^2$, we will sometimes write
\[E(\psi,\xi^2)=E(\psi),J(\psi,\xi^2)=J(\psi),\ \mathrm{and}\ -\lambda(\xi^2)=\inf_{H_0^1\cap H^2,\psi\neq0}\frac{E(\psi,\xi^2)}{J(\psi,\xi^2)}.\]

Hence, if $\mu_c>0$ and $\mu\in(0,\mu_c)$, for any $\xi^2\in[0,\xi_c^2)$, it holds that $\lambda(\xi^2)>0$.
\begin{rem}\label{r3.3}
{\it One should also notice that, under the assumptions that $\mu_c>0$ and $\mu\in(0,\mu_c)$, for any $\xi^2\geq\xi_c^2$, the fact
\[\xi^2\geq\xi_c^2=\sup_{H_0^1\cap H^2,\psi\neq0}\frac{\mathcal{N}_1(\psi)}{\mathcal{N}_2(\psi,\xi_c^2)}\geq\sup_{H_0^1\cap H^2,\psi\neq0}\frac{\mathcal{N}_1(\psi)}{\mathcal{N}_2(\psi,\xi^2)}\]
leads to
\[\xi^2\geq\frac{k_1(\psi^\prime(1))^2+k_0(\psi^\prime(0))^2-\mu\int_0^1({\psi}^{\prime\prime})^2}
{\mu\int_0^1\left(2({\psi}^\prime)^2+\xi^2{\psi}^2\right)},~\forall\psi\in H_0^1\cap H^2,\]
which further implies that $\lambda(\xi^2)\leq 0$ with $\lambda(\xi^2)=0$ if and only if $\xi^2=\xi_c^2.$

In addition, if $\mu\geq\mu_c$, then $\lambda(\xi^2)\leq 0$ for any $\xi^2\in[0,+\infty)$, and, $\lambda(\xi^2)=0$ if and only if $\mu=\mu_c$ and $\xi^2=0.$

In fact, one can see that
\[\lambda(\xi^2)\leq\sup_{H_0^1\cap H^2,\psi\neq0}\frac{(\mu_c-\mu)\int_0^1(\psi^{\prime\prime})^2}{2J(\psi,\xi^2)}-\mu\xi^2.\]

Therefore for any $\mu>\mu_c$, we have $\lambda(\xi^2)\leq -\mu<0$ if $\xi^2\geq 1$. Moreover, for any $\mu>\mu_c$ and $0\leq\xi^2\leq 1$, one has
\[\lambda(\xi^2)\leq\frac{(\mu_c-\mu)\|\psi^{\prime\prime}\|^2_{L^2}}{\|\psi\|^2_{H^1}}\leq\mu_c-\mu<0,\]
where (\ref{2.8}) has been used.

In conclusion, one gets that $\lambda(\xi^2)<\mu_c-\mu<0$ for any $\xi^2\in[0,+\infty)$ provided that $\mu>\mu_c.$}
\end{rem}

{\it Step 2.} In this step, we show that $\lambda$ is a bounded, continuous, strictly decreasing function with respect to $\xi^2$ on $[0,+\infty)$.
\begin{prop}\label{p4.2}
For $\mu_c>0$ and $\mu\in(0,\mu_c)$, the function $\lambda:[0,+\infty)\rightarrow\mathbb{R}$ is continuous, strictly decreasing and satisfies
\begin{align}\label{4.14}
\Lambda:=\max_{\xi^2\in[0,+\infty)}\lambda(\xi^2)=\lambda(0)\leq C_0.
\end{align}
where the constant $C_0$ is defined in (\ref{3.3}), which is positive in this case.
\end{prop}
\proof
In view of Remark \ref{r3.2} and similarly to Proposition \ref{p4.1}, for any $\xi_1^2,\xi^2_2\in[0,+\infty)$, we denote
\[\lambda(\xi_1^2)=\frac{-E(\psi_{\xi_1^2},\xi_1^2)}{J(\psi_{\xi_1^2},\xi_1^2)},
~\lambda(\xi_2^2)=\frac{-E(\psi_{\xi_2^2},\xi_2^2)}{J(\psi_{\xi_2^2},\xi_2^2)}.\]

Notice that
\begin{align}\label{4.16}
\lambda(\xi_1^2)&=\frac{-E(\psi_{\xi_1^2},\xi_1^2)}{J(\psi_{\xi_1^2},\xi_1^2)}
=\frac{-E(\psi_{\xi_1^2},\xi_2^2)}{J(\psi_{\xi_1^2},\xi_2^2)}+\frac{-E(\psi_{\xi_1^2},\xi_1^2)}{J(\psi_{\xi_1^2},\xi_1^2)}
-\frac{-E(\psi_{\xi_1^2},\xi_2^2)}{J(\psi_{\xi_1^2},\xi_2^2)}\nonumber\\
&\leq \lambda(\xi_2^2)
+\frac{\mathcal{N}_1(\psi_{\xi_1^2})(\xi_2^2-\xi_1^2)\int_0^1(\psi_{\xi_1^2})^2}{J(\psi_{\xi_1^2},\xi_1^2)J(\psi_{\xi_1^2},\xi_2^2)}
+\mu(\xi_2^2-\xi_1^2)+\frac{\mu(\xi_2^2-\xi_1^2)\left(\int_0^1(\psi^\prime_{\xi_1^2})^2\right)^2}
{J(\psi_{\xi_1^2},\xi_1^2)J(\psi_{\xi_1^2},\xi_2^2)}\nonumber\\
&\leq \lambda(\xi_2^2)+(C_0+2\mu)|\xi^2_2-\xi^2_1|,
\end{align}
where the fact that $\mathcal{N}_1(\psi_{\xi_1^2})>0$ has been used. The continuity of $\Lambda(\xi^2)$ then follows.

For $\xi_1^2<\xi_2^2$, by the definition of supremum, one has that
\begin{align}\label{4.15}
\lambda(\xi_1^2)=\frac{-E(\psi_{\xi_1^2},\xi_1^2)}{J(\psi_{\xi_1^2},\xi_1^2)}&\geq \frac{-E(\psi_{\xi_2^2},\xi_1^2)}{J(\psi_{\xi_2^2},\xi_1^2)}
=\frac{\mathcal{N}_1(\psi_{\xi_2^2})}{J(\psi_{\xi_2^2},\xi_1^2)}-\mu\xi_1^2
-\frac{\mu\xi_1^2\int_0^1(\psi^\prime_{\xi_2^2})^2}{J(\psi_{\xi_2^2},\xi_1^2)}\nonumber\\
&>\frac{\mathcal{N}_1(\psi_{\xi_2^2})}{J(\psi_{\xi_2^2},\xi_2^2)}-\mu\xi_2^2
-\frac{\mu\xi_2^2\int_0^1(\psi^\prime_{\xi_2^2})^2}{J(\psi_{\xi_2^2},\xi_2^2)}=\lambda(\xi_2^2),
\end{align}
where one has used the fact that $\mathcal{N}_1(\psi_{\xi_2^2})>0.$ This yields the monotonicity of $\lambda(\xi^2)$.

Consequently, $\Lambda=\lambda(0)$. Moreover, by using the same technique as in (\ref{3.2}), one can obtain that $\lambda(0)\leq C_0.$
\endproof

{{\it Step 3.} In this  step, we construct some growing mode solutions to (\ref{bdy1})-(\ref{lin}) by using the results in {\it Step 1} and {\it Step 2}.
\begin{prop}\label{p3.3}
{\it Let $f\in C_c^\infty(0,\xi_c^2)$ be a real-valued function and the real-valued functions $\phi(\xi^2,y)$, $\psi(\xi^2,y),\pi(\xi^2,y),\lambda(\xi^2)$ are the solutions, constructed in Proposition \ref{p3.1} and Proposition \ref{p3.2}, to problem (\ref{1.14})-(\ref{1.17}), where $\xi_c^2$ is the so called critical frequency which is positive and defined in (\ref{4.12}). Define
\begin{align}
&u^1(x,y,t)=-\frac{1}{2\pi}\int_{\mathbb{R}}f(\xi^2)i\phi(\xi^2,y)e^{\lambda(\xi^2)t}e^{ix\xi}d\xi,\label{4.191}\\
&u^2(x,y,t)=\frac{1}{2\pi}\int_{\mathbb{R}}f(\xi^2)\psi(\xi^2,y)e^{\lambda(\xi^2)t}e^{ix\xi}d\xi,\label{4.192}\\
&q(x,y,t)=\frac{1}{2\pi}\int_{\mathbb{R}}f(\xi^2)\pi(\xi^2,y)e^{\lambda(\xi^2)t}e^{ix\xi}d\xi.\label{4.193}
\end{align}
Then $(\mathbf{u}=(u^1,u^2),q)$ is a solution to linearized problem (\ref{bdy1})-(\ref{lin}). Due to the smoothness of functions $\phi(y),\psi(y), \pi(y)$, we also have the estimates
\begin{align}\label{4.20}
\|{\mathbf u}(0)\|_{H^k}+\|q(0)\|_{H^k}\leq \tilde{C}_k\left(\int_{\mathbb{R}}(1+\xi^2)^{k}|f(\xi)|^2\right)^{1/2}<+\infty, k\in\mathbb{N},
\end{align}
where constant $\tilde{C}_k>0$ depending on $k_0,k_1,\mu$ and $k$.

 Moreover, for every $t>0$, the boundedness of $\lambda(\xi^2)$ over $(0,\xi_c^2)$ implies that the solution $({\mathbf u}(t),q(t))\in H^k$ and satisfies
\begin{align}
&e^{\lambda_ft}\|{\mathbf{u}}(0)\|_{H^k}\leq\|{\mathbf{u}}(t)\|_{H^k}\leq e^{\Lambda t}\|{\mathbf{u}}(0)\|_{H^k},\label{4.21}\\
&e^{\lambda_ft}\|q(0)\|_{H^k}\leq\|q(t)\|_{H^k}\leq e^{\Lambda t}\|q(0)\|_{H^k},\label{4.23}
\end{align}
where
\begin{equation}\label{4.24}
\lambda_f:=\inf_{\xi^2\in\mathrm{supp}(f)}\lambda(\xi^2)>0
\end{equation}
and $\Lambda$ is a positive number defined in (\ref{4.14}).}
\end{prop}

\proof
It follows from (\ref{1.14})-(\ref{1.17}), Proposition \ref{p3.1}, \ref{p3.2}, Remark \ref{r3.3} and Proposition \ref{p4.2} that the solution given in (\ref{4.191})-(\ref{4.193}) satisfies (\ref{4.20})-(\ref{4.23}), and (\ref{4.24}) holds. This verification is similar to the proof of Theorem 2.4 in \cite{10}, and thus is omitted.

This completes the proof of linear instability part of Theorem \ref{thm1}.}
\endproof

\section{The nonlinear instability}
\subsection{Global existence and nonlinear energy estimates}

In this subsection, we prove that the nonlinear perturbed problem (\ref{nonl})-(\ref{bdy3}) admits at least one global strong solution.

The proof of local existence and uniqueness of strong solution is similar to that in section 4 of \cite{wtk} (see also section 2 of \cite{wx}). Therefore, in order to get the global existence of strong solutions, it suffices to derive some global energy estimates. To this end, let $({\mathbf u},q)$ be a strong solution of the perturbed problem (\ref{nonl})-(\ref{bdy3}). In the sequel, for simplicity, $C$ will denote a generic positive constant, which may depend on $k_1,k_0$ and $\mu$, and $C(\alpha,\beta)$ denotes some constant also depending on parameters $\alpha$ and $\beta.$

Testing $(\ref{nonl})_1$ by ${\mathbf u}$, integrating by part over $\Omega$ and using $(\ref{nonl})_2,$ boundary conditions (\ref{bdy1})-(\ref{bdy3}), one has that
\begin{align}\label{5.1}
\frac{1}{2}\frac{d}{dt}\int|{\mathbf u}(t)|^2+\mu\int|\nabla {\mathbf u}(t)|^2
=\int_{\mathbb{R}}\left(k_1|u^1(x,1)|^2+k_0|u^1(x,0)|^2\right):=I_1.
\end{align}

Notice that
\begin{align}\label{5.2}
I_1&=\int_0^1\frac{d}{dy}\left[\int_{\mathbb{R}}((k_1+k_0)y-k_0)|u^1(x,y)|^2\right]\nonumber\\
&\leq \frac{\mu}{2}\int|\partial_yu^1(x,y)|^2+2C_0\int|u^1(x,y)|^2.
\end{align}

Substituting (\ref{5.2}) into (\ref{5.1}) yields
\begin{align}\label{5.3}
\frac{d}{dt}\int|{\mathbf u}(t)|^2+\mu\int|\nabla {\mathbf u}(t)|^2\leq 4C_0\int|{\mathbf u}(t)|^2,
\end{align}
which, together with Gronwall inequality, implies that for any fixed $T>0$
\begin{align}\label{5.41111}
\sup_{0\leq t\leq T}\|{\mathbf u}(t)\|^2_{L^2}+\int_0^T\mu\|\nabla{\mathbf u}(t)\|_{L^2}^2dt\leq e^{4C_0T}\|{\mathbf u}_0\|^2_{L^2}.
\end{align}

Similarly, one gets that
\begin{align}\label{5.8}
\frac{1}{2}\frac{d}{dt}\int|{\mathbf u}_t(t)|^2+\mu\int|\nabla {\mathbf u}_t(t)|^2
\leq&4C_0\|{\mathbf u}_t\|_{L^2}^2+\frac{\mu}{4}\|\nabla{\mathbf u}_t\|^2_{L^2}-\int{\mathbf u}_t\cdot\nabla{\mathbf u}\cdot{\mathbf u}_t\nonumber\\
\leq&\int|{\mathbf u}_t|^2|\nabla{\mathbf u}|+4C_0\|{\mathbf u}_t\|_{L^2}^2+\frac{\mu}{4}\|\nabla{\mathbf u}_t\|^2_{L^2}\nonumber\\
\leq& 4C_0\|{\mathbf u}_t\|^2_{L^2}+C\|\nabla{\mathbf u}\|_{L^2}\|{\mathbf u}_t\|^2_{L^4}
+\frac{\mu}{4}\|\nabla{\mathbf u}_t\|^2_{L^2}\nonumber\\
\leq& \frac{\mu}{2}\|\nabla{\mathbf u}_t\|^2_{L^2}+C_1(1+\|\nabla{\mathbf u}\|^2_{L^2})\|{\mathbf u}_t\|^2_{L^2},
\end{align}
which implies that
\begin{align}\label{5.9}
\frac{d}{dt}\int|{\mathbf u}_t(t)|^2+\mu\int|\nabla {\mathbf u}_t(t)|^2\leq 2C_1(1+\|\nabla{\mathbf u}\|^2_{L^2})\|{\mathbf u}_t\|^2_{L^2}.
\end{align}

Multiplying $(\ref{nonl})_1$ by ${\mathbf u}_t $, integrating in space and recalling div${\mathbf u}_t=0$, one has
\begin{align}\label{5.12}
\int|{\mathbf u}_t(t)|^2=\int(\mu\Delta{\mathbf u}\cdot{\mathbf u}_t-{\mathbf u}\cdot\nabla{\mathbf u}\cdot{\mathbf u}_t)\lesssim\int(|{\mathbf u}||\nabla{\mathbf u}|+|\Delta{\mathbf u}|)|{\mathbf u}_t|.
\end{align}

Using the Cauchy inequality, the H\"older inequality and the Sobolev embedding inequalities, we arrive at
\begin{align}\label{5.13}
\|{\mathbf u}_t(t)\|_{L^2}^2&\lesssim\|{\mathbf u}(t)\|^2_{L^4}\|\nabla{\mathbf u}(t)\|^2_{L^4}+\|\nabla^2{\mathbf u}(t)\|^2_{L^2}\nonumber\\
&\lesssim\|{\mathbf u}(t)\|^2_{L^2}\|{\mathbf u}(t)\|^2_{H^2}+\|{\mathbf u}(t)\|^2_{H^2}\nonumber\\
&\leq C(1+e^{4C_0t})\|{\mathbf u}(t)\|^2_{H^2}.
\end{align}

Taking $t\rightarrow 0^+$ in the above inequality yields
\begin{align}\label{5.14}
\limsup_{t\rightarrow 0^+}\|{\mathbf u}_t(t)\|_{L^2}^2\leq  C\|{\mathbf u}_0\|^2_{H^2},
\end{align}
where $C>0$ depends also on $\|{\mathbf u}_0\|^2_{L^2}$.

Therefore, applying Gronwall inequality to (\ref{5.9}), we have
\begin{align}\label{5.91}
\sup_{0\leq t\leq T}\|{\mathbf u}_t(t)\|^2_{L^2}+\int_0^T\mu\|{\nabla\mathbf u}_t(t)\|_{L^2}^2dt\leq C(T,\|{\mathbf u}_0\|_{H^2}).
\end{align}

Testing $(\ref{nonl})_1$ by ${\mathbf u}_t $, integrating by part over $\Omega$ and using $(\ref{nonl})_2,$ boundary conditions (\ref{bdy1})-(\ref{bdy3}), we obtain
\begin{align}\label{5.4}
\frac{\mu}{2}\frac{d}{dt}\int|\nabla{\mathbf u}(t)|^2+\int|{\mathbf u}_t(t)|^2=
&\int_{\mathbb{R}}\left(k_1u^1(x,1)u_t^1(x,1)+k_0u^1(x,0)u_t^1(x,0)\right)\nonumber\\
&-\int{\mathbf u}\cdot\nabla{\mathbf u}\cdot{\mathbf u}_t:=I_2+I_3.
\end{align}

Similar to (\ref{5.2}), it holds that
\begin{align}\label{5.5}
I_2&\leq C\int(|{\mathbf u}||{\mathbf u}_t|+|\nabla{\mathbf u}||{\mathbf u}_t|+|{\mathbf u}||\nabla{\mathbf u}_t|)
\leq C(\|{\mathbf u}_t\|^2_{L^2}+\|{\mathbf u}\|^2_{H^1}+\|\nabla{\mathbf u}_t\|^2_{L^2}),
\end{align}
and
\begin{align}\label{5.6}
I_3=\int{\mathbf u}\otimes{\mathbf u}:\nabla{\mathbf u}_t\leq\|{\mathbf u}\|^2_{L^4}\|\nabla{\mathbf u}_t\|_{L^2}\leq C(\|{\mathbf u}\|^2_{L^2}\|\nabla{\mathbf u}_t\|^2_{L^2}+\|{\mathbf u}\|^2_{H^1}),
\end{align}
where the two-dimensional interpolation inequality
$$\|{\mathbf u}\|^2_{L^4}\lesssim \|{\mathbf u}\|_{L^2}\|{\mathbf u}\|_{H^1}$$ has been used.

Substituting (\ref{5.5}) and (\ref{5.6}) into (\ref{5.4}), and integrating over $[0,t]\subset[0,T]$, one gets
\begin{align}\label{5.7}
\sup_{0\leq t\leq T}\|\nabla{\mathbf u}(t)\|^2_{L^2}+\int_0^T\mu\|{\mathbf u}_t(t)\|_{L^2}^2dt\leq C(T,\|{\mathbf u}_0\|_{H^2}).
\end{align}

Finally, we recall that the pair $({\mathbf u},q)$ solves the Stokes equations
\begin{align}\label{5.16}
\begin{cases}
-\mu\Delta{\mathbf u}+\nabla q=-{\mathbf u}_t-{\mathbf u}\cdot\nabla{\mathbf u},\\
\mathrm{div}{\mathbf u}=0.
\end{cases}
\end{align}

By Stokes estimate (\ref{A2}) in the Appendix, it is clear that
\begin{align}\label{5.161}
\|\nabla^2{\mathbf u}\|_{L^2}^2+\|\nabla q\|_{L^2}^2&\leq C\|{\mathbf u}_t+{\mathbf u}\cdot\nabla{\mathbf u}\|^2_{L^2}+C\|{\mathbf u}\|^2_{L^2}\nonumber\\
&\leq C\|{\mathbf u}_t\|^2_{L^2}+C\|{\mathbf u}\|^2_{L^4}\|\nabla{\mathbf u}\|^2_{L^4}+C\|{\mathbf u}\|^2_{L^2}\nonumber\\
&\leq C\left(\|{\mathbf u}_t\|^2_{L^2}+\|{\mathbf u}\|^2_{L^2}\right)+C\|{\mathbf u}\|^3_{H^1}\|{\mathbf u}\|_{H^2}\nonumber\\
&\leq C(T,\|{\mathbf u}_0\|_{H^2})+\frac{1}{2}\|\nabla^2{\mathbf u}\|_{L^2}^2,
\end{align}
which implies
\begin{align}\label{5.17}
\|\nabla^2{\mathbf u}\|_{L^2}^2+\|\nabla q\|_{L^2}^2&\leq C(T,\|{\mathbf u}_0\|_{H^2}).
\end{align}

Summing up, we have obtained the global energy estimates to guarantee the global existence of strong solutions (see Proposition \ref{p5.1}) as follows:
\begin{align}\label{5.180}
&\sup_{0\leq t\leq T}(\|{\mathbf u}(t)\|_{H^2}^2+\|\nabla q(t)\|_{L^2}^2+\|{\mathbf u}_t(t)\|_{L^2}^2)
+\int_0^T(\|\nabla{\mathbf u}(t)\|_{L^2}^2+\|{\mathbf u}_t(t)\|_{H^1}^2)dt\nonumber\\
\leq& C(T,\|{\mathbf u}_0\|_{H^2}).
\end{align}
\begin{prop}\label{p5.1}
For any given $T>0$ and initial data ${\mathbf u}_0\in H^2$ satisfying the compatibility condition $\mathrm{div}{\mathbf u}_0=0$, there exists a strong solution $({\mathbf u},q)\in C([0,T];H^2\times H^1)$ to the perturbed problem (\ref{nonl})-(\ref{bdy3}). Moreover, there exists a constant $\bar{\sigma}\in (0,1]$, such that
\begin{align}\label{5.20}
&\|{\mathbf u}(t)\|^2_{H^2}+\|({\mathbf u}_t,\nabla q)(t)\|^2_{L^2}+\int_0^t\|(\nabla{\mathbf u},{\mathbf u}_t,\nabla{\mathbf u}_t)(s)\|^2_{L^2}ds\nonumber\\
\leq& C_1\left(\|{\mathbf u}_0\|^2_{H^2}+\int_0^t\|{\mathbf u}(s)\|^2_{L^2}ds\right),
\end{align}
provided that $\|{\mathbf u}(t)\|^2_{H^2}\leq\bar{\sigma}$ on $[0,T]$. Here the constant $C_1$ depends only on $k_0,k_1$, and $\mu.$
\end{prop}
\proof
One can follow the proof of section 4 of \cite{wtk} (see also section 2 of \cite{wx}) to get the local existence and uniqueness of strong solution to the nonlinear perturbed problem (\ref{nonl})-(\ref{bdy3}). Then the global existence and uniqueness of the strong solution can be shown easily by using the above global {\it a priori} estimate (\ref{5.180}).

It remains to prove (\ref{5.20}). In view of the assumption that $\|{\mathbf u}(t)\|^2_{H^2}\leq\bar{\sigma}$, one can estimate $I_2,I_3$ in (\ref{5.4}) as follows.
\begin{align}
I_2&\leq C\int(|{\mathbf u}||{\mathbf u}_t|+|\nabla{\mathbf u}||{\mathbf u}_t|+|{\mathbf u}||\nabla{\mathbf u}_t|)
\leq \frac{1}{4}\|{\mathbf u}_t\|^2_{L^2}+C_{\epsilon}\|{\mathbf u}\|^2_{H^1}+\epsilon\|\nabla{\mathbf u}_t\|^2_{L^2},\label{5.200}\\
I_3&\leq \int|{\mathbf u}||\nabla{\mathbf u}||{\mathbf u}_t|\leq\|{\mathbf u}_t\|_{L^2}\|\nabla{\mathbf u}\|_{L^4}\|{\mathbf u}\|_{L^4}\leq \frac{1}{4}\|{\mathbf u}_t\|^2_{L^2}+C\|{\mathbf u}\|^2_{H^2}\|{\mathbf u}\|^2_{L^2}.\label{5.201}
\end{align}

Substituting (\ref{5.200}) and (\ref{5.201}) into (\ref{5.4}),we have
\begin{align}\label{5.202}
\frac{d}{dt}\int\mu|\nabla{\mathbf u}(t)|^2+\int|{\mathbf u}_t|^2\leq C_{\epsilon}\|{\mathbf u}\|^2_{H^1}+\epsilon\|\nabla{\mathbf u}_t\|^2_{L^2}+C\|{\mathbf u}\|^2_{H^2}\|{\mathbf u}\|^2_{L^2}.
\end{align}

Adding $K_1\times(\ref{5.202}),(\ref{5.9}),$ and $K_2\times(\ref{5.3})$ up with suitable large $K_1>0,K_2>0$ and taking $\epsilon>0$ small enough, we arrive at
\begin{align}\label{5.10}
&\frac{d}{dt}\|\left({\mathbf u},\sqrt{\mu}\nabla{\mathbf u},{\mathbf u}_t\right)(t)\|^2_{L^2}
+\|\left(\sqrt{\mu}\nabla{\mathbf u},{\mathbf u}_t,\sqrt{\mu}\nabla{\mathbf u}_t\right)\|^2_{L^2}\nonumber\\
\leq&C\|{\mathbf u}\|^2_{L^2}+C\|{\mathbf u}\|^2_{H^2}\left(\|{\mathbf u}\|^2_{L^2}
+\|{\mathbf u}_t\|^2_{L^2}\right)\nonumber\\
\leq&C\|{\mathbf u}\|^2_{L^2}+C\bar{\sigma}\left(\|{\mathbf u}\|^2_{L^2}
+\|{\mathbf u}_t\|^2_{L^2}\right)
\end{align}
provided that $\|{\mathbf u}\|^2_{H^2}\leq\bar{\sigma}$.

Then, for suitably small $\bar{\sigma}\in(0,1]$, one can get that
\begin{align}\label{5.11}
\frac{d}{dt}\|\left({\mathbf u},\sqrt{\mu}\nabla{\mathbf u},{\mathbf u}_t\right)(t)\|^2_{L^2}
+\|\left(\sqrt{\mu}\nabla{\mathbf u},{\mathbf u}_t,\sqrt{\mu}\nabla{\mathbf u}_t\right)\|^2_{L^2}\leq C\|{\mathbf u}\|^2_{L^2}.
\end{align}

Thus, it follows from (\ref{5.11}) that
\begin{align}\label{5.15}
\|\left({\mathbf u},\nabla{\mathbf u},{\mathbf u}_t\right)(t)\|^2_{L^2}
+\int_0^t\|\left(\nabla{\mathbf u},{\mathbf u}_t,\nabla{\mathbf u}_t\right)(s)\|^2_{L^2}ds\leq C\left(\|{\mathbf u}_0\|^2_{H^2}+\int_0^t\|{\mathbf u}(s)\|^2_{L^2}ds\right)
\end{align}

Moreover, under the assumption $\|{\mathbf u}\|^2_{H^2}\leq\bar{\sigma}\leq 1$ as in the proof of (\ref{5.161}), one can get that
\begin{align}\label{5.18}
\|\nabla^2{\mathbf u}\|_{L^2}^2+\|\nabla q\|_{L^2}^2\leq C(\|{\mathbf u}_t\|^2_{L^2}+\|{\mathbf u}\cdot\nabla{\mathbf u}\|^2_{L^2}+{\|{\mathbf u}\|^2_{L^2}})\leq C\|({\mathbf u},\nabla{\mathbf u},{\mathbf u}_t)\|^2_{L^2},
\end{align}
which, together with (\ref{5.15}), implies (\ref{5.20}).
\endProof

\subsection{Proof of Theorem \ref{thm2} (i): nonlinear instability}
In this subsection, we apply the bootstrap argument proposed by Y. Guo et al. in [2] to prove the nonlinear instability. More precisely, we shall show that there exists a constant $\varepsilon>0$ such that for any $\delta>0$, there exists a solution $\mathbf{u}^\delta(t)$ to the nonlinear problem (\ref{nonl})-(\ref{bdy3}) with initial data $\|\mathbf{u}^\delta_0\|_{H^2}=\delta$ and an escape time $T^\delta>0$ such that $\|\mathbf{u}^\delta(T^\delta)\|_{H^2}>\varepsilon.$

To this end, we first give the following elementary inequality, which will be used in this section and in the next section.
\begin{prop}\label{p5.2}
Let $\mathbf{w}\in H_\sigma^1(\Omega)\cap H^2(\Omega)$, then it holds that
\begin{align}\label{Key}
-\mu\int|\nabla \mathbf{w}|^2dxdy+k_1\int_{\mathbb{R}}|w^1(x,1)|^2dx+k_0\int_{\mathbb{R}}|w^1(x,0)|^2dx\leq\Lambda\int|\mathbf{w}|^2dxdy,
\end{align}
where $\Lambda$ is defined in (\ref{4.14}).
\end{prop}
\proof
For any function $g\in L^2(\Omega)$, let
 \[\hat{g}(\xi,y)=\int_{\mathbb{R}}g(x,y)e^{-i\xi x}dx,\ \xi\in\mathbb{R}.\]
 Then it follows from {\it Fubini} theorem and {\it Parseval} equality that $\hat{g}\in L^2(\Omega)$ and
 \begin{align}\label{5.350}
 \int_{\Omega}|g(x,y)|^2dxdy=\frac{1}{2\pi}\int_{\Omega}|\hat{g}(\xi,y)|^2d\xi dy.
 \end{align}
Hence,
\begin{align}\label{5.351}
&-2\pi\int\mu|\nabla{\mathbf w}|^2dxdy+2\pi\int_{\mathbb{R}}k_1|w^1(x,1)|^2dx+2\pi\int_{\mathbb{R}} k_0|w^1(x,0)|^2dx\nonumber\\
=&-\int_{\mathbb{R}}\int_0^1\mu\left[|i\xi\hat{w}^1(\xi,y)|^2+|i\xi\hat{w}^2(\xi,y)|^2
+|\partial_y\hat{w}^1(\xi,y)|^2+|\partial_y\hat{w}^2(\xi,y)|^2\right]d\xi dy\nonumber\\
&+\int_{\mathbb{R}}\left[k_1|\hat{w}^1(\xi,1)|^2+k_0|\hat{w}^1(\xi,0)|^2\right]d\xi.
\end{align}

 For simplicity, denoting $\phi(y)=i\hat{w}^1(\xi,y),\psi(y)=\hat{w}^2(\xi,y)$ for fixed $\xi\neq 0$, then (\ref{5.351}) becomes
\begin{align}\label{5.352}
&-2\pi\int\mu|\nabla{\mathbf w}|^2dxdy+2\pi\int_{\mathbb{R}}k_1|w^1(x,1)|^2dx+2\pi\int_{\mathbb{R}} k_0|w^1(x,0)|^2dx\nonumber\\
=&-\int_{\mathbb{R}}\mu\int_0^1(|\xi\phi|^2+|\xi\psi|^2+|\phi^\prime|^2+|\psi^\prime|^2)dy-(k_1|\phi(1)|^2+k_0|\phi(0)|^2)d\xi,
\end{align}
where $^\prime=\partial_y.$
Set
\begin{align*}
Z(\phi,\psi;\xi)=-\mu\int_0^1(|\xi\phi|^2+|\xi\psi|^2+|\phi^\prime|^2+|\psi^\prime|^2)dy+k_1|\phi(1)|^2+k_0|\phi(0)|^2.
\end{align*}

Clearly,
\[Z(\phi,\psi;\xi)=Z(\mathfrak{R}\phi,\mathfrak{R}\psi;\xi)+Z(\mathfrak{I}\phi,\mathfrak{I}\psi;\xi).\]
 Thus, it suffice to bound $Z$ when $\phi,\psi$ are real-value functions.

  Notice that $\mathrm{div}\mathbf{w}=0$, so $\xi\phi+\psi^\prime=0.$ Then, using (\ref{2.2}), we may rewrite
  \[Z(\phi,\psi;\xi)=-2E(\psi,\xi)/\xi^2,\ \xi\neq 0\]
 and hence it follows from the definition and Proposition \ref{p4.2} that
\begin{align}\label{5.353}
Z(\phi,\psi;\xi)\leq \frac{2\lambda(\xi^2)}{\xi^2}J(\psi;\xi^2)=\frac{\lambda(\xi^2)}{\xi^2}\int_0^1(\xi^2|\psi|^2+|\psi^\prime|^2)dy
\leq\Lambda\int_0^1(|\psi|^2+|\phi|^2)dy.
\end{align}

Translating this inequality back to the original form yields that
\begin{align}\label{5.354}
\mu\int_0^1(|i\xi\hat{w}^1(\xi,y)|^2+|i\xi\hat{w}^2(\xi,y)|^2+|\partial_y\hat{w}^1(\xi,y)|^2+|\partial_y\hat{w}^2(\xi,y)|^2)dy
\nonumber\\
+k_1|\hat{w}^1(\xi,1)|^2+k_0|\hat{w}^1(\xi,0)|^2\leq\Lambda\int_0^1(|i\hat{w}^1|^2+|\hat{w}^2|^2)dy.
\end{align}

Then, integrating each side of this inequality over all $\xi\in\mathbb{R}$ and using (\ref{5.350}), we obtain (\ref{Key}). The Proposition follows.
\endProof

Now we are on the position to prove the nonlinear instability.

By Theorem \ref{thm1}, one can construct a solution to the linear problem (\ref{bdy1})-(\ref{lin}) in the form:
\begin{equation}\label{5.21}
\bar{{\mathbf u}}(x,y,t)=
\begin{cases}
\bar{u}^1(x,y,t)=-\frac{1}{2\pi}\int_{\mathbb{R}}f(\xi^2)i\phi(\xi^2,y)e^{\lambda(\xi^2)t}e^{ix\xi}d\xi\\
\bar{u}^2(x,y,t)=\frac{1}{2\pi}\int_{\mathbb{R}}f(\xi^2)\psi(\xi^2,y)e^{\lambda(\xi^2)t}e^{ix\xi}d\xi
\end{cases}
\in H^2
\end{equation}
with initial data
\begin{equation}\label{5.210}
\bar{\mathbf u}_0(x,y)=
\begin{cases}
\bar{u}^1(x,y,0)=-\frac{1}{2\pi}\int_{\mathbb{R}}f(\xi^2)i\phi(\xi^2,y)e^{ix\xi}d\xi\\
\bar{u}^2(x,y,0)=\frac{1}{2\pi}\int_{\mathbb{R}}f(\xi^2)\psi(\xi^2,y)e^{ix\xi}d\xi
\end{cases}
\in H^2
\end{equation}
 satisfying $\mathrm{div}\bar{\mathbf u}_0=0$ and $\|\bar{\mathbf u}_0\|_{H^2}=1$.

  Moreover, one can suitably choose the cut-off function $f\in C_0^\infty(0,\xi_c^2)$ such that
 \begin{equation}\label{5.211}
 \lambda_*\leq\lambda_f<\Lambda,
 \end{equation}
 where $\lambda_f$ and $\Lambda$ are defined in (\ref{4.24}) and (\ref{4.14}), and $\lambda_*>\frac{\Lambda}{2}$ will be determined later.

Denote ${\mathbf u}_0^\delta:=\delta\bar{\mathbf{u}}_0$ and $C_2:=\|\bar{\mathbf{u}}_0\|_{L^2}$. By {Proposition \ref{p5.1}}, for any $\delta\in(0,\bar{\sigma})$, there exists a global strong solution $({\mathbf u}^\delta,p^\delta)\in C([0,T]; H^2\times H^1)$ to (\ref{nonl})-(\ref{bdy3}), with the initial data ${\mathbf{u}}_0^\delta$ satisfying $\|{\mathbf{u}}_0^\delta\|_{H^2}=\delta.$

Then, for any $\delta\in(0,\bar{\sigma})$ such that $\delta<\varepsilon_0$, define
\begin{align}\label{5.22}
T^\delta:=\frac{1}{\lambda_*}\ln{\frac{\varepsilon_0}{\delta}}~\mathrm{i.e.}~\delta e^{\lambda_* T^\delta}=\varepsilon_0,
\end{align}
where $\varepsilon_0>0$, independent of $\delta$, is a small constant to be determined, and $\lambda_*=\lambda_*(\varepsilon_0,\delta)$ is the same parameter as in (\ref{5.211}).

Furthermore, define
\begin{align}\label{5.23}
T^*=\sup\{t\in(0,+\infty)\big|\|\mathbf{u}^\delta\|_{H^2}\leq\bar{\sigma}\}
\end{align}
and
\begin{align}\label{5.24}
T^{**}=\sup\{t\in(0,+\infty)\big|\|\mathbf{u}^\delta\|_{L^2}\leq 2C_2\delta e^{\lambda_* t}\}.
\end{align}

Obviously, $T^*,T^{**}>0$ and
\begin{align}
&\|{\mathbf u}^{\delta}(T^*)\|_{H^2}=\bar{\sigma},~\mathrm{if}~ T^*<+\infty,\label{5.25}\\
&\|{\mathbf u}^{\delta}(T^{**})\|_{L^2}=2C_2\delta e^{\lambda_* T^{**}},~\mathrm{if}~ T^{**}<+\infty.\label{5.26}
\end{align}

For any $t\leq\min\{T^*,T^{**},T^\delta\}$, (\ref{5.20}) implies that
\begin{align}\label{5.27}
\|\mathbf{u}^\delta(t)\|^2_{H^2}+\|\mathbf{u}_t^\delta(t)\|^2_{L^2}&\leq C_1\|{\mathbf u}_0^\delta\|^2_{H^2}+C_1\int_0^t(2C_2\delta e^{\lambda_* s})^2ds\nonumber\\
&\leq C_1\delta^2+2C_1C_2^2\delta^2e^{2\lambda_* t}/\lambda_*:= C_3\delta^2e^{2\lambda_* t},
\end{align}
where $C_3$, independent of $\delta$, is a positive constant.

Denote $\mathbf{u}^d=\mathbf{u}^\delta-\delta\bar{\mathbf{u}}$ and $\mathbf{u}_\delta^L=\delta\bar{\mathbf{u}}$. Note that $\mathbf{u}_\delta^L$ is also a strong solution to the linearized problem (\ref{bdy1})-(\ref{lin}) with the initial data $\mathbf{u}_0^\delta\in H^2$. Thus $\mathbf{u}^d$ solves
\begin{align}\label{5.28}
\begin{cases}
\mathbf{u}^d_t+\nabla p^d-\mu\Delta\mathbf{u}^d=-\mathbf{u}^\delta\cdot\nabla\mathbf{u}^\delta,\\
\mathrm{div}\mathbf{u}^d=0,
\end{cases}
\end{align}
with the boundary conditions
\begin{align}
&u^{d,2}(x,0)=u^{d,2}(x,1)=0,\ x\in \mathbb{R},\label{5.29}\\
&\partial_yu^{d,1}(x,1)=\frac{k_1}{\mu}u^{d,1}(x,1),\ x\in \mathbb{R},\label{5.30}\\
&\partial_yu^{d,1}(x,0)=-\frac{k_0}{\mu}u^{d,1}(x,0),\ x\in \mathbb{R},\label{5.31}
\end{align}
where $u^{d,1}$ and $u^{d,2}$ stand for the first and second component of $\mathbf{u}^d$ respectively, and the initial condition $\mathbf{u}^d(0)=\mathbf{0}$.

Multiplying $(\ref{5.28})_1$ by ${\mathbf u}^d$ gives that
\begin{align}\label{5.32}
\frac{1}{2}\frac{d}{dt}\int|\mathbf{u}^d|^2=-\mu\int|\nabla{\mathbf u}^d|^2+\sum^1_{i=0}\int_{\mathbb{R}}k_i|u^{d,1}(x,i)|^2
-\int{\mathbf u}^\delta\cdot\nabla{\mathbf u}^\delta\cdot{\mathbf u}^d.
\end{align}

Notice that
\begin{align}\label{5.33}
\int{\mathbf u}^\delta\cdot\nabla{\mathbf u}^\delta\cdot{\mathbf u}^d
\leq\int|{\mathbf u}^\delta\cdot\nabla{\mathbf u}^\delta||{\mathbf u}^d|
\leq\|{\mathbf u}^\delta\cdot\nabla{\mathbf u}^\delta\|_{L^2}\|{\mathbf u}^d\|_{L^2}
\leq C_4\|{\mathbf u}^\delta\|^2_{H^2}\|{\mathbf u}^d\|_{L^2}.
\end{align}

In addition, Proposition \ref{p5.2} implies that
\begin{align}\label{5.35}
-\mu\int|\nabla{\mathbf u}^d|^2+\sum^1_{i=0}\int_{\mathbb{R}}k_i|u^{d,1}(x,i)|^2\leq\Lambda\int|{\mathbf u}^d|^2,
\end{align}
where $\Lambda>0$ is defined in (\ref{4.14}).

Substituting $(\ref{5.33})$ and $(\ref{5.35})$ into (\ref{5.32}) gives that
\begin{align}\label{5.36}
\frac{d}{dt}\|\mathbf{u}^d\|_{L^2}\leq \Lambda\|{\mathbf u}^d\|_{L^2}+C_4\|{\mathbf u}^\delta\|^2_{H^2}.
\end{align}

Thus, it follows from the Gronwall inequality, (\ref{5.27}) and (\ref{5.36}) that
\begin{align}\label{5.37}
\|{\mathbf u}^d\|_{L^2}\leq C_4e^{\Lambda t}\int_0^te^{-\Lambda s}\|{\mathbf u}^\delta(s)\|^2_{H^2}ds\leq C_3C_4\delta^2e^{\Lambda t}\int_0^te^{(2\lambda_*-\Lambda)s}ds\leq C_5\delta^2e^{2\lambda_* t},
\end{align}
where the condition $2\lambda_*-\Lambda>0$ has been used.

Now we claim that
\begin{align}\label{5.38}
T^\delta=\min\{T^\delta,T^*,T^{**}\},\ \mathrm{provided}\ \varepsilon_0=\min\{\frac{\bar{\sigma}}{2\sqrt{C_3}},\frac{C_2}{4C_5}\}.
\end{align}

Indeed, if $T^*=\min\{T^\delta,T^*,T^{**}\}$, then $T*<+\infty$. It follows from (\ref{5.27}) and (\ref{5.22}) that
\begin{align}\label{5.39}
\|{\mathbf u}^\delta(T^*)\|_{H^2}\leq\sqrt{C_3}\delta e^{\lambda_* T^*}\leq\sqrt{C_3}\delta e^{\lambda_* T^\delta}=\sqrt{C_3}\varepsilon_0<\bar{\sigma},
\end{align}
which contradicts (\ref{5.25}).

If $T^{**}=\min\{T^\delta,T^*,T^{**}\}$, then $T^{**}<+\infty$. In view of (\ref{4.21}),(\ref{5.22}) and (\ref{5.37}), one obtains that
\begin{align}\label{5.40}
\|{\mathbf u}^\delta(T^{**})\|_{L^2}&\leq\|{\mathbf u}^L_\delta(T^{**})\|_{L^2}+\|{\mathbf u}^d(T^{**})\|_{L^2}
\leq C_2\delta e^{\Lambda T^{**}}+C_5\delta^2 e^{2\lambda_*{T^{**}}}\nonumber\\
&\leq C_2\delta e^{\lambda_* T^{**}}\left(e^{(\Lambda-\lambda_*)T^\delta}+\frac{C_5}{C_2}\delta e^{\lambda_*{T^\delta}}\right)
 \leq C_2\delta e^{\lambda_* T^{**}}\left[\left(\frac{\varepsilon_0}{\delta}\right)^{\frac{\Lambda}{\lambda_*}-1}+\frac{1}{4}\right].
\end{align}

Take
\begin{align}\label{5.401}
\lambda_*=\Lambda\ln{\left(\frac{2\varepsilon_0}{\delta}\right)}/\ln{\left(\frac{5\varepsilon_0}{2\delta}\right)}.
\end{align}

Then $\lambda_*>\Lambda/2$ since that $\varepsilon_0>\delta$. Therefore,
\begin{align}\label{5.400}
\|{\mathbf u}^\delta(T^{**})\|_{L^2}<2C_2\delta e^{\lambda_* T^{**}},
\end{align}
which contradicts (\ref{5.26}). Therefore, (\ref{5.38}) holds.

Finally, we use (\ref{4.21}),(\ref{5.22}) and (\ref{5.37}) to deduce that
\begin{align}\label{5.41}
\|{\mathbf u}^\delta(T^\delta)\|_{L^2}&\geq\|{\mathbf u}^L_\delta(T^\delta)\|_{L^2}-\|{\mathbf u}^d(T^\delta)\|_{L^2}
\geq C_2\delta e^{\lambda_* T^\delta}-C_5\delta^2e^{2\lambda_* T^\delta}>C_2\varepsilon_0/2,
\end{align}
which completes the proof of {Theorem \ref{thm2} (i)} by defining $\varepsilon:=C_2\varepsilon_0/2.$

\section{The linear and nonlinear stability}
{In the first subsection, we will prove Theorem \ref{thm3}, namely, asymptotic stability of the linear and nonlinear system under the assumption of $\mu>\mu_c\geq 0$. We will analyse for the case $\mu\geq\mu_c\geq 0$ in the second subsection to complete the proof of the stability part of Theorem \ref{thm1} and the proof of Theorem \ref{thm2} (ii).}

It follows from Proposition \ref{p5.2} and Remark \ref{r3.3} that for any $\mathbf{u}(t)\in H_\sigma^1\cap H^2$, it holds that
\begin{align}\label{6.1}
\int_{\mathbb{R}}\left[k_1|u^1(x,1))|^2+k_0|u^1(x,0)|^2\right]
-\mu\int|\nabla{\mathbf u}|^2\leq\Lambda\int|{\mathbf u}|^2,
\end{align}
where $\Lambda< 0$ provided $\mu>\mu_c$, while $\Lambda=0$ for $\mu=\mu_c>0$. This is crucial for the proof of the stability. In what follows, for simplicity, we denote by $C$ a generic positive constant, which may depend on $k_1,k_0$ and $\mu$.

\subsection{Proof of Theorem \ref{thm3}.}

\hskip0.7cm{\it Proof of Theorem \ref{thm3} (i): general initial data.}

Standard energy estimates and (\ref{6.1}) yield
\begin{align}\label{6.3}
\frac{1}{2}\frac{d}{dt}\|\mathbf{u}\|_{L^2}^2=\sum^1_{i=0}\int_{\mathbb{R}}k_i|u^1(x,i)|^2-\mu\|\nabla\mathbf{u}\|_{L^2}^2\leq \Lambda\|\mathbf{u}\|_{L^2}^2,
\end{align}
where $\Lambda<0.$ This implies that
\begin{align}\label{6.02}
\|\mathbf{u}(t)\|_{L^2}\leq e^{\Lambda t}\|\mathbf{u}_0\|_{L^2}.
\end{align}

In addition, one has
\begin{align}\label{6.03}
\frac{1}{2}\frac{d}{dt}\|\mathbf{u}\|_{L^2}^2+(\mu-\mu_c)\|\nabla\mathbf{u}\|_{L^2}^2=\sum^1_{i=0}\int_{\mathbb{R}}k_i|u^1(x,i)|^2
-\mu_c\|\nabla\mathbf{u}\|_{L^2}^2\leq 0,
\end{align}
which gives that
\begin{align}\label{6.04}
\|\mathbf{u}(t)\|^2_{L^2}+\int_0^t\|\nabla\mathbf{u}(t)\|^2_{L^2}\leq C\|\mathbf{u}_0\|^2_{L^2}.
\end{align}

Applying $\partial_t$ to $(\ref{nonl})_1$, taking the inner product of the result with ${\mathbf u}_t$, and treating the boundary terms as in (\ref{6.03}), one gets that for any $\epsilon>0$,
\begin{align}\label{6.05}
&\frac{1}{2}\frac{d}{dt}\|\mathbf{u}_t\|_{L^2}^2+(\mu-\mu_c)\|\nabla\mathbf{u}_t\|_{L^2}^2\lesssim \int|{\mathbf u}_t|^2|\nabla{\mathbf u}|\nonumber\\
\lesssim&\|\mathbf{u}_t\|^2_{L^4}\|\nabla\mathbf{u}\|_{L^2}
\lesssim \|\mathbf{u}_t\|_{L^2}\|\mathbf{u}_t\|_{H^1}\|\nabla\mathbf{u}\|_{L^2}\leq \epsilon\|\mathbf{u}_t\|^2_{H^1}
+C_\epsilon\|\nabla\mathbf{u}\|^2_{L^2}\|\mathbf{u}_t\|^2_{L^2},
\end{align}
where H\"older inequality, Young inequality and Sobolev embedding theorems have been used.

It follows from (\ref{5.4}) and a similar argument as for (\ref{6.05}) that
\begin{align}\label{6.06}
\frac{\mu}{2}\frac{d}{dt}\|\nabla\mathbf{u}\|_{L^2}^2+\|\mathbf{u}_t\|^2_{L^2}
&\lesssim\int(|\mathbf{u}||\nabla\mathbf{u}||\mathbf{u}_t|+|\mathbf{u}_t||\mathbf{u}|
+|\nabla\mathbf{u}_t||\mathbf{u}|+|\mathbf{u}_t||\nabla\mathbf{u}|)\nonumber\\
&\lesssim \|\mathbf{u}\|_{L^4}\|\nabla\mathbf{u}\|_{L^2}\|\mathbf{u}_t\|_{L^4}+\|\mathbf{u}\|_{H^1}\|\mathbf{u}_t\|_{H^1}\nonumber\\
&\leq \epsilon\|\mathbf{u}_t\|^2_{H^1}+C_\epsilon(\|\mathbf{u}\|^4_{H^1}+\|\mathbf{u}\|^2_{H^1}).
\end{align}

Adding (\ref{6.05}) and (\ref{6.06}) and taking $\epsilon$ small enough yield
\begin{align}\label{6.07}
\frac{d}{dt}\left(\|\nabla\mathbf{u}\|^2_{L^2}+\|\mathbf{u}_t\|^2_{L^2}\right)+\|\mathbf{u}_t\|^2_{H^1}\leq C\|\nabla\mathbf{u}\|^2_{L^2}\left(\|\nabla\mathbf{u}\|^2_{L^2}+\|\mathbf{u}_t\|^2_{L^2}\right)+C\|\mathbf{u}\|^2_{H^1}.
\end{align}

Notice that
\begin{align}\label{6.08}
\int_0^t\|\mathbf{u}(s)\|^2_{H^1}\leq\int_0^te^{\Lambda s}\|\mathbf{u}_0\|^2_{L^2}ds+\int_0^t\|\nabla\mathbf{u}(s)\|^2_{L^2}ds\leq C\|\mathbf{u}_0\|^2_{L^2},
\end{align}
where (\ref{6.02}) and (\ref{6.04}) have been used.

Thus, by applying Gronwall inequality to (\ref{6.07}), one gets
\begin{align}\label{6.09}
\|\nabla\mathbf{u}(t)\|^2_{L^2}+\|\mathbf{u}_t(t)\|^2_{L^2}+\int_0^t\|\mathbf{u}_t(s)\|^2_{H^1}ds\leq C\|\mathbf{u}_0\|^2_{H^2},
\end{align}
where (\ref{5.14}) has been used.

By the Stokes estimate (\ref{A3}), we have
\begin{align}\label{6.010}
\|\mathbf{u}(t)\|_{H^2}+\|\nabla q(t)\|_{L^2}&\lesssim \|\mathbf{u}_t(t)\|_{L^2}+\|\mathbf{u}(t)\cdot\nabla\mathbf{u}(t)\|_{L^2}\nonumber\\
&\lesssim\|\mathbf{u}_t(t)\|_{L^2}+\|\mathbf{u}(t)\|_{L^\infty}\|\nabla\mathbf{u}(t)\|_{L^2}\nonumber\\
&\lesssim\|\mathbf{u}_t(t)\|_{L^2}+\|\mathbf{u}(t)\|^{1/2}_{L^2}\|\mathbf{u}(t)\|^{1/2}_{H^2}\|\nabla\mathbf{u}(t)\|_{L^2}\nonumber\\
&\leq\frac{1}{2}\|\mathbf{u}(t)\|_{H^2}+C\|\mathbf{u}_t(t)\|_{L^2}+C\|\mathbf{u}(t)\|_{L^2}\|\nabla\mathbf{u}(t)\|^2_{L^2},
\end{align}
which, together with (\ref{6.04}) and (\ref{6.09}), implies that
\begin{align}\label{6.011}
\|\mathbf{u}(t)\|_{H^2}\leq C\|\mathbf{u}_0\|_{H^2}.
\end{align}

Furthermore, interpolation inequality implies that
\begin{align}\label{6.012}
\|\nabla\mathbf{u}(t)\|_{L^2}\leq C\|\mathbf{u}(t)\|^{1/2}_{H^2}\|\mathbf{u}(t)\|^{1/2}_{L^2}\leq Ce^{\Lambda t/2},
\end{align}
which, together with (\ref{6.02}), yields the third inequality of (\ref{6.01}) by taking $\alpha=-\Lambda/2$.

Notice that (\ref{6.09}) implies $\|\mathbf{u}_t(s)\|^2_{L^2}\rightarrow 0$ as $t\rightarrow+\infty$. Then one can see from (\ref{6.010}) and (\ref{6.012}) that
$\|\mathbf{u}(t)\|_{H^2}\rightarrow 0$ as $t\rightarrow+\infty.$

Theorem \ref{thm3} (i) follows.

{\it Proof of Theorem \ref{thm3} (ii): small initial data.}

In fact, replacing $\mathbf{u}$ in (\ref{6.1}) by $\mathbf{u}_t$ and integrating by part, one can  re-estimate (\ref{6.05}) as
\begin{align}\label{6.013}
&\frac{1}{2}\frac{d}{dt}\|\mathbf{u}_t\|_{L^2}^2+(-\Lambda)\|\mathbf{u}_t\|^2_{L^2}\leq\int|{\mathbf u}_t||\nabla{\mathbf u}_t||\mathbf{u}|\nonumber\\
&\lesssim\|\nabla\mathbf{u}_t\|_{L^2}\|{\mathbf u}_t\|_{L^2}\|{\mathbf u}\|_{L^\infty}
\lesssim\|\nabla\mathbf{u}_t\|_{L^2}\|{\mathbf u}_t\|_{L^2}\|{\mathbf u}\|_{H^2}\nonumber\\
&\leq \epsilon\|\nabla\mathbf{u}_t\|^2_{L^2}+C_\epsilon\|{\mathbf u}_t\|_{L^2}^2\|{\mathbf u}\|^2_{H^2}.
\end{align}

Adding (\ref{6.013}) and (\ref{6.05}) with $\epsilon>0$ small enough, we obtain that
\begin{align}\label{6.014}
\frac{d}{dt}\|\mathbf{u}_t\|_{L^2}^2+\frac{(-3\Lambda)}{4}\|\mathbf{u}_t\|^2_{L^2}\leq C\|{\mathbf u}_t\|_{L^2}^2\|{\mathbf u}\|^2_{H^2}.
\end{align}

Taking $\|{\mathbf u}_0\|_{H^2}$ sufficiently small in (\ref{6.011}), we have $C\|{\mathbf u}(t)\|^2_{H^2}\leq-\Lambda/4$. Then, we have
\[\frac{d}{dt}\|\mathbf{u}_t(t)\|_{L^2}^2\leq \Lambda/2\|\mathbf{u}_t(t)\|^2_{L^2}\]
which implies that
\begin{align}\label{6.015}
\|\mathbf{u}_t(t)\|_{L^2}\leq Ce^{\Lambda t/2}.
\end{align}
(\ref{6.2}) follows by (\ref{6.015}),(\ref{6.010}) and (\ref{6.012}), where $\beta=-\Lambda/2.$

Theorem \ref{thm3} (ii) is proved.
\endProof

\subsection{Proof of the stability part of Theorem \ref{thm1} and Theorem \ref{thm2}(ii)}

In this subsection, we prove the stability part of Theorem \ref{thm1}. For Theorem \ref{thm2}(ii), we will give a remark at the end.

{\it Proof of the stability part of Theorem \ref{thm1}.}

{\it Step 1.} $\mu>\mu_c$: decay estimates.

In fact, one can see that in the linearized situation, by similar energy method as used in the proof of Theorem \ref{thm3}, it is easy to obtain the decay rate that
\[\|\mathbf{u}(t)\|_{H^2}^2\leq Ce^{\Lambda t/2}\|\mathbf{u}_0\|_{H^2}^2,\]
which automatically implies that
$\|\mathbf{u}(t)\|_{H^2}^2\rightarrow 0\ \mathrm{as}\ t\rightarrow 0,$
since $\Lambda<0$ provided $\mu>\mu_c$. It should be noticed that in the linearized situation, the initial data need not to be small for us to obtain this decay estimate.

{\it Step 2.} $\mu=\mu_c$: continuous dependence on initial data.

Similarly to (\ref{6.3}), since in this case, $\Lambda=0$, one only has
\begin{align}\label{6.9}
\frac{d}{dt}\|\mathbf{u}\|_{L^2}^2=\int_{\mathbb{R}}\left[k_1|u^1(x,1)|^2+k_0|u^1(x,0)|^2\right]-\mu_c\|\nabla\mathbf{u}\|_{L^2}^2\leq 0.
\end{align}

Multiplying $(\ref{lin})_1$ by $\mathbf{u}_t$, using $(\ref{lin})_2$ and the boundary conditions yield
\begin{align}\label{6.10}
\frac{1}{2}\frac{d}{dt}\left(\mu_c\|\nabla\mathbf{u}\|_{L^2}^2-\sum_{i=0}^1k_i\int_{\mathbb{R}}|u^1(x,i)|^2\right)
+\|\mathbf{u}_t\|_{L^2}^2=0.
\end{align}

Similar to (\ref{5.9}) in Section 4, one has
\begin{align}\label{6.11}
\frac{d}{dt}\|{\mathbf u}_t(t)\|_{L^2}^2+\mu\|\nabla {\mathbf u}_t(t)\|_{L^2}^2\leq 2C_0\|{\mathbf u}_t\|^2_{L^2}.
\end{align}

Adding up $K_4\times(\ref{6.10}),K_5\times(\ref{6.9})$ and (\ref{6.11}) with suitably large $K_4>0$, we arrive at
\begin{align}\label{6.15}
\frac{d}{dt}\left(\|(\sqrt{K_5}\mathbf{u},\sqrt{\mu_cK_4}\nabla\mathbf{u},\mathbf{u}_t)(t)\|^2_{L^2}
-K_4\sum_{i=0}^1k_i\int_{\mathbb{R}}|u^1(x,i,t)|^2\right)\leq 0.
\end{align}

Finally, integrating (\ref{6.15}) over $(0,t)$, we obtain
\begin{align}\label{6.16}
\|(\sqrt{K_5}\mathbf{u},\sqrt{\mu_cK_4}\nabla\mathbf{u},\mathbf{u}_t)(t)\|^2_{L^2}
-K_4\sum_{i=0}^1k_i\int_{\mathbb{R}}|u^1(x,i,t)|^2\leq C\|\mathbf{u}_0\|^2_{H^2}.
\end{align}

Taking $K_5>0$ large enough and applying the Stokes estimates (\ref{A2}) imply that
\begin{align}\label{6.17}
\|\mathbf{u}(t)\|^2_{H^2}+\int_0^t\|\mathbf{u}_t(s)\|^2_{H^1}ds\leq C\|\mathbf{u}_0\|^2_{H^2}.
\end{align}

{The stability part of Theorem \ref{thm1} is proved.}
\endproof

\begin{rem}\label{r5.1}
 {\it In this remark, we state a proof of Theorem \ref{thm2} (ii). In fact, for the general case  $\mu\geq\mu_c$, it follows from (\ref{6.02}) that $\|\mathbf{u}(t)\|_{L^2}^2\leq \|\mathbf{u}_0\|_{L^2}^2$, which completes the proof of Theorem \ref{thm2} (ii).}
\end{rem}

\section*{Appendix}
\appendix
\section{The Stokes estimates}
{Denote that $\Gamma_i:=\mathbb{R}\times\{i\},i=0,1$, and $\Omega:=\mathbb{R}\times(0,1)$. Consider the following Stokes equations with Navier-slip boundary conditions,
\begin{equation}\label{A1}
\begin{cases}
-\mu\Delta\mathbf{u}-\nabla p=\mathbf{F},&\Omega,\\
\mathrm{div}\mathbf{u}=0,&\Omega,\\
u^2=0,&\Gamma_0\cup\Gamma_1,\\
\mu\partial_yu^1-k_1u^1=0,&\Gamma_1,\\
\mu\partial_yu^1+k_0u^1=0,&\Gamma_0.\\
\end{cases}
\end{equation}
\begin{thm}\label{thm4}
Suppose that $\mu>0$, $\mathbf{F}\in H^{m-1},\ m\in\mathbb{N}$ and $(\mathbf{u},p)$ solves (\ref{A1}). Then the following claims holds.

(i) If $\mu>\mu_c>0$, then
\begin{equation}\label{A3}
\|\mathbf{u}\|^2_{H^{m+1}}+\|\nabla p\|^2_{H^{m-1}}\leq C\|\mathbf{F}\|^2_{H^{m-1}},
\end{equation}
where $\mu_c$ is defined in (\ref{2.4}) and $C$ is a positive constant depending only on $\mu,k_1,k_0,m.$

(ii) If $0<\mu\leq\mu_c$, it holds that
\begin{equation}\label{A2}
\|\mathbf{u}\|^2_{H^{m+1}}+\|\nabla p\|^2_{H^{m-1}}\leq C\left(\|\mathbf{F}\|^2_{H^{m-1}}+\|\mathbf{u}\|^2_{L^2}\right),
\end{equation}
where $C$ is a positive constant depending only on $\mu,k_1,k_0,m.$
\end{thm}

{\it Proof of (i) ($\mu>\mu_c>0$).}

{\it Step 1.} Multiplying $(\ref{A1})_1$ by $\mathbf{u}$ and integrating by part over $\Omega$, we have
\begin{align*}
\mu\int|\nabla\mathbf{u}|^2-k_1\int_{\Gamma_1}|u^1|^2-k_0\int_{\Gamma_0}|u^1|^2=\int\mathbf{F}\cdot\mathbf{u}.
\end{align*}

Since $\mu>\mu_c$, one may choose $\delta>0$ such that $\mu-\mu_c-\delta>0$. Rewriting above equality as
\begin{align*}
(\mu-\mu_c-\delta)\int|\nabla\mathbf{u}|^2+(\mu_c+\delta)\int|\nabla\mathbf{u}|^2-k_1\int_{\Gamma_1}|u^1|^2-k_0\int_{\Gamma_0}|u^1|^2=\int\mathbf{F}\cdot\mathbf{u},
\end{align*}
and using Proposition \ref{p5.2}, one gets
\begin{align*}
(\mu-\mu_c-\delta)\int|\nabla\mathbf{u}|^2-\Lambda\int|\mathbf{u}|^2\leq\int\mathbf{F}\cdot\mathbf{u}
\end{align*}
where $\Lambda<0$ is a constant which implies
\begin{align}\label{A4}
\|\mathbf{u}\|_{H^1}\leq C\|\mathbf{F}\|_{H^{-1}}.
\end{align}

{\it Step 2.} Applying horizontal differential operator $\nabla^m_x$ to (\ref{A1}), one has, similar to (\ref{A4}), that
\begin{align}\label{A5}
\|\nabla^m_x\mathbf{u}\|_{H^1}&\leq C\|\nabla_x^m\mathbf{F}\|_{H^{-1}}\leq C\|\mathbf{F}\|_{H^{m-1}}.
\end{align}

{\it Step 3.} Since that $\partial\Omega$ is horizontally flat, one has
\begin{align*}
\|\mathbf{u}\|^2_{H^m(\partial\Omega)}=\|\mathbf{u}\|^2_{L^2(\partial\Omega)}+\|\nabla^m_x\mathbf{u}\|^2_{L^2(\partial\Omega)}.
\end{align*}

Then, the trace theorem and (\ref{A4})-(\ref{A5}) yield
\begin{align}\label{A6}
\|\mathbf{u}\|^2_{H^{m+\frac{1}{2}}(\partial\Omega)}&=\|\mathbf{u}\|^2_{H^{\frac{1}{2}}(\partial\Omega)}
+\|\nabla^m_x\mathbf{u}\|^2_{H^{\frac{1}{2}}(\partial\Omega)}\nonumber\\
&\leq C\left(\|\mathbf{u}\|^2_{H^1}+\|\nabla^m_x\mathbf{u}\|^2_{H^1}\right)\nonumber\\
&\leq C\|\mathbf{F}\|^2_{H^{m-1}}.
\end{align}

{\it Step 4.} By the regularity of $\mathbf{u}$ on the boundary, (\ref{A6}), one may use the classical estimates for the following problem
\begin{equation}\label{A7}
\begin{cases}
-\mu\Delta\mathbf{u}-\nabla p=\mathbf{F},&\Omega,\\
\mathrm{div}\mathbf{u}=0,&\Omega,\\
\mathbf{u}=\mathbf{u},&\partial\Omega
\end{cases}
\end{equation}
to obtain the following inequality (see \cite{Temam,Ladyz})
\begin{align}\label{A8}
\|\mathbf{u}\|_{H^{m+1}}+\|\nabla p\|_{H^{m-1}}&\leq C\left(\|\mathbf{F}\|_{H^{m-1}}+\|\mathbf{u}\|_{H^{m+\frac{1}{2}}(\partial\Omega)}\right)
\leq C\|\mathbf{F}\|_{H^{m-1}}
\end{align}
where (\ref{A6}) has been used. Thus (i) is proved.

{\it Proof of (ii) ($0<\mu\leq\mu_c$).}

 In this case, we also have
\begin{align*}
\mu\int|\nabla\mathbf{u}|^2-k_1\int_{\Gamma_1}|u^1|^2-k_0\int_{\Gamma_0}|u^1|^2=\int\mathbf{F}\cdot\mathbf{u}
\end{align*}
which implies that
 \begin{align*}
\mu\int|\nabla\mathbf{u}|^2&\leq|k_1|\int_{\Gamma_1}|u^1|^2+|k_0|\int_{\Gamma_0}|u^1|^2+\int\mathbf{F}\cdot\mathbf{u}\\
&=\int_{\mathbb{R}}dx\int_0^1\left[\left((|k_0|+|k_1|)y-|k_0|\right)(u^1)^2\right]^{\prime}_ydy+\int\mathbf{F}\cdot\mathbf{u}\\
&\leq \frac{\mu}2\int|\nabla\mathbf{u}|^2+C\int|\mathbf{u}|^2+C\|\mathbf{F}\|_{H^{-1}(\Omega)}^2.
\end{align*}

 Therefore,
\begin{align}\label{A9}
\|\mathbf{u}\|_{H^1}\leq C\left(\|\mathbf{F}\|_{H^{-1}}+\|\mathbf{u}\|_{L^{2}}\right).
\end{align}

Now claim (ii) follows from (\ref{A9}) and the similar steps in the proof of (i).

Theorem \ref{A1} follows.
\endProof

\noindent{\bf Final Remark.~}This is an erratum of the previous version with the same title (arXiv:1608.03019), which is published in Journal of Mathematical Fluid Mechanics, 20(2018),603-629. Since the value of $\mu_c$ in the case of $k_0=k_1$ should be revised, we refer the readers to Remark \ref{r2.3} in this version for the details. The erratum for Section 2 has been submitted to Journal of Mathematical Fluid Mechanics.

\acknowledgment

The authors would like to thank Professor Yanjin Wang and Professor Huanyao Wen for their constructive comments and discussions in the preparation of the paper.

The authors would also like to thank Dr. Tien-Tai Nguyen for pointing out the mistake of the value of $\mu_c$ at $k_1=k_0>0$ in the previous version(arXiv:1608.03019) of this paper in \cite{2}.

Ding's reseach is supported by the National Natural Science Foundation of China (No.11071086, No.11371152, No.11128102 and No.11571117).

 Xin's research is partially supported by Zheng Ge Ru Foundation, Hong Kong RGC Earmarked Research Grants CUHK-14305315 and CUHK4048/13P, NSFC/RGC Joint Research Scheme Grant N-CUHK 443-14, and a Focus Area Grant from the Chinese University of Hong Kong.


\begin{thebibliography}{99}


\bibitem{1} Y. Achdou, O. Pironneau, F. Valentin, Effective boundary conditions for laminar flow over periodic rough boundaries,
 {\it J. Comput. Phys.} {\bf 147}(1998), 187-218.

\bibitem{2} C. Amrouche, A. Rejaiba, $L^p$-theory for Stokes and Navier-Stokes equations with Navier boundary condition, {\it J. Diff. Equs}, {\bf 256}, 2014, 1515-1547.
\bibitem{Amrouche} C. Amrouche, N. H. Seloula, On the Stokes equations with the Navier-type boundary conditions, {\it Diff. Equns. Appl.}, {\bf 3}(4), 2011, 581-607.
\bibitem{3} S. Antontsev, H. de Oliveira, Navier-Stokes equations with absorption under slip boundary conditions: existence,
uniqueness and extinction in time, {\it RIMS K\^{o}ky\^{u}roku Bessatsu} {\bf B1}(2007) 21-41.

\bibitem{4} E. B\"ansch, Finite element discretization of the Navier-Stokes equations with free capillary surface, {\it Numer. Math.} {\bf 88}, 2001, 203-235.

\bibitem{6} G. Beavers, D. Joseph, Boundary conditions at a naturally permeable wall, {\it J. Fluid Mech.} {\bf 30}, 1967, 197-207.

\bibitem{7} S. Chandrasekhar, {\it Hydrodynamic and Hydromagnetic Stability}, Internat. Ser. Monogr. Phys., Clarendon Press, Oxford, 1961.
\bibitem{1993}D. S. Chauhan, K. S. Shekhawat, Heat transfer in Couette flow of a compressible Newtonian fluid in the presence of a naturally permeable boundary, {\it J. Phys.D: Appl. Phys.} {\bf 26}, 1993, 933-936.

\bibitem{8} P. Drazin, W. Reid, {\it Hydrodynamic Stability}, Section Edition, Cambridge University Press, 2004.
\bibitem{gie} G. M. Gie, J. P. Kelliher, Boundary layer analysis of the Navier-Stokes equations with generalized Navier boundary conditions, {\it J. Diff. Equns.}, {\bf 253}, 2012, 1862-1892.
\bibitem{9} Y. Guo, C. Hallstrom, D. Spirn, Dynamics Near Unstable, Interfacial Fluids, {\it Commun. Math. Phys.} {\bf 270}, 2007, 635-689.
\bibitem{12} Y. Guo, Y. Han, Critical Rayleigh Number in Rayleigh-B\'enard Convection, {\it Quart. of Appl. Math.}, {\bf 68}(1), 2010, 149-160.
\bibitem{10} Y. Guo, I. Tice, Linear Rayleigh-Taylor Instability for Viscous, Compressible Fluid, {\it SIMA J. Math.Anal.} {\bf 42}(4), 2010, 1688-1720.

\bibitem{11} Y. Guo, I. Tice, Compressible, inviscid Rayleigh-Taylor instability, {\it Indiana Univ. Math. J.} {\bf 60}(2), 2011, 677-712.


\bibitem{guoyan2017} Y. Guo, I. Tice, Stability of contact line in fluids: 2D STOKES flow, 2017, arXiv:1603.03721v1.




\bibitem{haase} A. S. Haase, J. A. Wood, R. G. H. Lammertink, J. H. Snoeijer, Why bumpy is better: the role of the dissipaption distribution in slip flow over a bubble mattress, {\it Phys. Rev. Fluid}, {\bf 1}, 2016, 054101.


\bibitem{13} W. J\"ager, A. Mikeli\'c, On the Roughness-induced effective boundary conditions for an incompressible viscous flow,
{\it J. Diff. equs.} {\bf170}, 2001, 96-122.

\bibitem{14} W. J\"ager, A. Mikeli\'c, On the interface boundary condition of Beavers, Joseph, and Saffman,
 {\it SIAM J. Appl. Math.} {\bf 60}, 2000, 1111-1127.

\bibitem{15} F. Jiang, S. Jiang, G. Ni, Nonlinear instability for nonhomogeneous incompressible viscous fluids, {\it Sci. China Math.}, {\bf 56}(4), 2013, 665-686.

\bibitem{16} F. Jiang, S. Jiang, On Instability and Stability of Three-dimensional Gravity Driven Viscous Flows in a Boundary Domain,
{\it Adv. Math.}, {\bf 264}, 2014, 831-863.

\bibitem{17} V. John, Slip with friction and penetration with resistance boundary conditions for the Navier-Stokes equation-numerical test and aspect of the implementation, {\it J. Comput. Appl. Math.},  147, 2002, 287-300.

\bibitem{kelliher} J. P. Kelliher, Navier-Stokes equations with Navier boundary conditions for a bounded domain in plane, {\it SIAM J. Math. Anal.}, {\bf 38}(1), 2006, 210-232.
\bibitem{Ladyz} O. A. Ladyzhenskaya, Mathematical Theory of Viscous Incompressible Flow, Gordon, 1969.


\bibitem{lizhang} H. Li, X. Zhang, Stability of plane Couette flow for the compressible Navier-Stokes equations with Navier-slip boundary, 2016, Preprint.

\bibitem{1995} J. Magnaudet, M. Riverot, J. Fabre, Accelerated flows past a rigid sphere or a spherical bubble. Part 1. Steady straining flow, {\it J. Fluid Mech.} {\bf 284}, 1995, 97-135.

\bibitem{18} C. Navier, Sur les lois de l\'equilibre et du mouvement des corps \'elastiques, {\it Mem. Acad. R. Sci. Inst. France}
     {\bf 6}, 1827, 369.

\bibitem{19} T. Qian, X. Wang, P. Sheng, Molecular scale contact line hydrodynamics of immiscible flows, {\it Physical Review E} {\bf 68}, 2003, 016306.
\bibitem{TTN}Tien-Tai Nguyen, Linear and Nonlinear Analysis of the Rayleigh-Taylor system with Navier-slip Boundary conditions, arXiv:2204.09857v1(2022)
\bibitem{191} J. Serrin, Mathematical Principles of Classical Fluid Mechanics, {\it Encyclopedia of
Physics} VIII/1, Springer-Verlag, Berlin, 1959, 125-263).

\bibitem{20} V. Solonnikov, V. \v{S}\v{c}adilov, A certain boundary value problem for the stationary system of Navier-Stokes equations, {\it Trudy Mat. Inst. Steklov.}, 125, 1973, 196-210; translation in {\it Proc. Steklov Inst. Math.}, {\bf 125}, 1973, 186-199.


\bibitem{Temam}    R. Temam, Navier-Stokes equations, Studies in Mathematics and its applications
2, North-Holland, Amsterdam, 1984.

\bibitem{21} H. B. da Veiga, On the regularity of flows with Ladyzhenskaya shear-dependent viscosity and slip or nonslip boundary conditions, {\it Comm. Pure Appl. Math.}, {\bf LVIII}, 2005, 552-577.

\bibitem{22} Y. Wang, I. Tice, The Viscous Surface-Internal Wave Problem: Nonlinear Rayleigh-Taylor Instability, {\it Comm. P.D.E.}, {\bf 37}, 2012, 1967-2028.




\bibitem{wtk} Yanjin Wang, Ian Tice, Chanwoo Kim, The viscous surface-internal wave problem: global well-posedness and decay, Arch. Rational Mech. Anal., {\bf 212}, 2014, 1-92.

\bibitem{wx} Yanjin Wang, Zhouping Xin, Vanishing viscosity and surface tension limits of incompressible viscous surface waves, arXiv:1504.00152.






\bibitem{23} Y. Xiao, Z. Xin, On the Vanishing Viscosity Limit for the 3D Navier-Stokes Equations with a Slip Boundary Condition, {\it Comm. Pure Appl. Math.}, {\bf 60}, 2007, 1027-1055.

\bibitem{24} Y. Xiao, Z. Xin, On the Inviscid Limit of the 3D Navier-Stokes Equations with Generalized Navier-slip Boundary Conditions, {\it Comm. Math. Stat.}, {\bf 1}(3), 2013, 259-279.



\end{thebibliography}
\end{document}